\documentclass{article}

\usepackage[english]{babel}

\usepackage[letterpaper,top=2cm,bottom=2cm,left=3cm,right=3cm,marginparwidth=1.75cm]{geometry}

\usepackage{amsmath}
\usepackage{graphicx}
\usepackage[colorlinks=true, allcolors=blue]{hyperref}

\usepackage{authblk}
\usepackage{lipsum}
\usepackage[utf8]{inputenc}
\usepackage{graphicx}
\usepackage{epstopdf}
\usepackage{algorithmic}
\ifpdf
  \DeclareGraphicsExtensions{.eps,.pdf,.png,.jpg}
\else
  \DeclareGraphicsExtensions{.eps}
\fi
\usepackage[caption=false]{subfig}
\usepackage{amsmath,amssymb,amsfonts}
\usepackage{verbatim}
\usepackage{xcolor}
\usepackage{multirow}
\usepackage[square,numbers]{natbib}

\newcommand{\dff}{\partial}

\newcommand{\BS}{\boldsymbol}
\newcommand{\bs}{\boldsymbol}
\newcommand{\rmd}{\mathrm{d}}
\newcommand{\rmi}{\mathrm{i}}
\newcommand{\rminc}{\mathrm{inc}}
\newcommand{\rmsc}{\mathrm{sc}}
\newcommand{\rmtot}{\mathrm{tot}}
\newcommand{\rmtr}{\mathrm{tr}}
\newcommand{\rmin}{\mathrm{in}}
\newcommand{\rmout}{\mathrm{out}}

\title{Helmholtz equation and non-singular boundary elements applied to multi-disciplinary physical problems}
\author[1]{Evert Klaseboer}
\author[2, *]{Qiang Sun}
\affil[1]{Institute of High Performance Computing, 1 Fusionopolis Way, Singapore 138632, Singapore}
\affil[2]{Australian Research Council Centre of Excellence for Nanoscale BioPhotonics, School of Science, RMIT University, Melbourne, VIC 3001, Australia}

\affil[*]{qiang.sun@rmit.edu.au}

\date{}

\begin{document}
\maketitle

\begin{abstract}
The famous scientist Hermann von Helmholtz was born 200 years ago. Many complex physical wave phenomena in engineering can effectively be described using one or a set of equations named after him: the Helmholtz equation. Although this has been known for a long time from a theoretical point of view, the actual numerical implementation has often been hindered by divergence free and/or curl free constraints. There is further a need for a numerical method that is accurate, reliable and takes into account radiation conditions at infinity. The classical boundary element method (BEM) satisfies the last condition, yet one has to deal with singularities in the implementation. We review here how a recently developed singularity-free three-dimensional (3D) boundary element framework with superior accuracy can be used to tackle such problems only using one or a few Helmholtz equations with higher order (quadratic) elements which can tackle complex curved shapes. Examples are given for acoustics (a Helmholtz resonator among others) and electromagnetic scattering. 
\end{abstract}


\section{Introduction}

Hermann von Helmholtz (1821-1894)~\cite{HelmholtzBook,HelmholtzBook2} was born 200 years ago on 31 August 1821. The scalar Helmholtz equation, $\nabla^2 \phi + k^2 \phi = 0$, which inherits its name from the famous German scientist, appears in many fields of science and engineering involving waves. The waves thus described can be longitudinal (as in sound waves), transverse (as in electromagnetic waves) or both (as in dynamic elastic waves in solid materials). This appears to us as the perfect moment to revisit the Helmholtz equation and see how it can be solved most efficiently for real three dimensional engineering problems.

The main goal of this article is to demonstrate that many classical engineering problems such as in acoustics, electromagnetics and elasticity can, in principle, all be formulated using one or a set of scalar Helmholtz equations. Thus a robust and efficient method to solve the Helmholtz equation is most desirable. Having developed a boundary element method for the scalar Helmholtz equation which is free of singularities, it is now possible to tackle more problems in classical applied physics essentially using the same numerical framework. This method was developed after the realization that, ideally, if there are no singularities in the physical system, no such singularities should have to appear in the mathematical equivalent description of this physical system. We will show examples mainly for (decaying) scattering waves from an object situated in an infinite medium and focus on sound and electromagnetic waves. Especially for the field of electromagnetics, the framework presented here is substantially different from the industry standard implementation. 

The structure of this paper proceeds as follows. After this short introduction, different classic physical problems that can be formulated based on the scalar Helmholtz equation are demonstrated in Sec.~\ref{sec:ScalarHelmholtz}, with as examples, among others, the Helmholtz resonator in acoustics. In Sec.~\ref{sec:BEMChapter}, the recently developed non-singular boundary integral method to solve the Helmholtz equation is derived in detail, in which all the singularities in the integrands and the terms associated with the solid angle are eliminated analytically. After a brief discussion on curl free vector fields in Sec.~\ref{sec:Curlfree}, we discuss divergence free vector fields in Sec.~\ref{sec:divfree}, with particular attention to electromagnetic scattering problems. Then, the discussion and conclusions are given in Secs.~\ref{sec:Discussion} and~\ref{sec:Conclusion}, respectively.

\section{Scalar Helmholtz equation}\label{sec:ScalarHelmholtz}
\subsection{Acoustics}

The wave equation describing the propagation of a quantity $\phi'$ is in the time domain $
\nabla^2 \phi' = 
\frac{1}{c^2}
\frac{\partial^2 \phi'}{\partial t^2}$,
with $c$ the wave speed and $t$ time. Assuming a harmonic time dependence $\phi' =\phi \exp({-\rmi\omega t})$, with angular frequency $\omega$ and unit imaginary number `$\rmi$', the wave equation transforms into the well-known Helmholtz equation as:
\begin{align}
\nabla^2 \phi + k^2 \phi = 0 \label{eq:Helm}
\end{align}
with $k=\omega/c$. This is the equation as it appeared in Helmholtz' book~\cite{HelmholtzBook} on ``acoustics in pipes with open ends'' while studying organ pipes as Eq. (3b) on page 18/19. The Helmholtz equation is often used to describe acoustic waves, in which $\phi$ represents the velocity potential or the pressure perturbation in the medium. Normally, $k$ is a real valued number,
often referred to as wave number in wave physics, while $\phi$ is a complex quantity.
The boundary conditions are usually that $\phi$ is given (Dirichlet condition), or its normal derivative $\rmd \phi/\rmd n= \boldsymbol n \cdot \nabla \phi$ is given (Neumann condition) with $\boldsymbol n$ the normal vector pointing out of the domain, or a combination of those two (Robin conditions). Sound waves are typical examples of so-called longitudinal waves. 

\subsection{Non-singular boundary element method for the Helmholtz equation}\label{sec:BEMChapter}

Only for a very limited number of cases an exact theoretical solution can be found for the Helmholtz equation. Therefore an efficient numerical method to solve it is required. A good candidate is the boundary integral method (or boundary element method when it is used in discretized form) as it can represent surface geometry accurately and reduces a three-dimensional (3D) problem to a two-dimensional (2D) problem~\cite{Banerjee1981}.

The Helmholtz equation is an elliptic equation. This means that the properties on the surface(s) $S$ of the domain determine entirely the solution anywhere in the domain. This allows us to write the following surface integral equation which is equivalent to Eq.~(\ref{eq:Helm}):
\begin{equation}
    c(\BS{x}_0) \phi(\BS{x}_{0}) + \int_{S} \phi(\BS{x}) \frac{\dff{G_{k}(\BS{x},\BS{x}_{0})}}{\dff{n}} \text{ d} S(\BS{x}) = \int_{S}  \frac{\dff{\phi(\BS{x})}}{\dff{n}} G_{k}(\BS{x},\BS{x}_{0}) \text{ d} S(\BS{x}) \label{eq:BEM}
\end{equation}
where $G_k \equiv G_k(\BS x,\BS x_0) = \exp{(\rmi k |\BS{x}-\BS{x}_{0}|)}/|\BS{x}-\BS{x}_{0}|$ is the Green's function, $\boldsymbol x_0$ is a source point on the surface, $\boldsymbol x$ is the computation point on the surface $S$, and $c(\BS{x}_0)$ is a constant related to the solid angle at $\boldsymbol x_0$.

Though the conventional boundary integral method (CBIM) in Eq.~(\ref{eq:BEM}) is often used to solve the Helmholtz equation, Eq.~(\ref{eq:Helm}), there are two main problems worth mentioning. Firstly, the integrands of CBIM in Eq.~(\ref{eq:BEM}) are singular when $\BS{x}\rightarrow\BS{x}_0$. Although these singular integrands can be integrated, their practical treatment in numerical implementations is not straightforward~\cite{Gaul2003, Vijayakumar1989, Rosen1992, Rosen1995}. Such singularities have no physical basis but are a purely mathematical artefact. Secondly, the computation of the solid angle $c(\BS{x}_0)$ is also not straightforward for collocation nodes, which makes the use of surface elements other than planar constant elements rather tedious. 

Treatment and even elimination of the singularities has been a subject of intense study from the beginning of the boundary element method. Well known are the so-called constant potential subtraction methods that remove the singularity on the left hand side of Eq.~(\ref{eq:BEM}) (see for example \cite{Seybert1985}), yet the singular behavior on the right hand side of Eq.~(\ref{eq:BEM}) still has to be dealt with. Another way to avoid the singularity issue is to use the method of fundamental solutions or equivalent source method~\cite{Fairweather2003,Lee2017} where the sources are put outside of the domain. 

We have developed an advanced non-singular boundary element method in which the solid angle and singularities in the integrals are removed analytically. This improvement has made the boundary element method much easier to implement, especially for higher order quadratic elements. 

Recently, the non-singular boundary integral methods (NS-BIM) or boundary regularised integral equation formulations (BRIEF) have been developed for applications in fluid mechanics~\cite{KlaseboerJFM2012, SunPRE2013, SunEABE2014, SunPoF2015}, acoustics~\cite{SunRoySoc2015, KlaseboerJASA2017, 2104.13137}, molecular electrostatics~\cite{SunJCP2016}, electromagnetics~\cite{KlaseboerIEEETAP2017, SunPRB2017, KlaseboerAO2017}, interactions between light and matter~\cite{Sun2021, Sun2022}  and linear elastic waves~\cite{KlaseboerJE2018}. The objective of the non-singular boundary integral method is to analytically remove the singularities in $G_{k}(\BS{x},\BS{x}_{0})$ and $\dff{G_{k}}(\BS{x},\BS{x}_{0})/\dff{n}$ as $\BS{x}\rightarrow \BS{x}_{0}$ as well as the solid angle $c(\BS{x}_0)$. 

The singular behavior of $G_{k}(\BS{x},\BS{x}_{0})$ is the same as that of the free-space Green's function of the Laplace equation: $ G_{0} \equiv G_{0}(\BS{x},\BS{x}_{0}) = {1}/{|\BS{x}-\BS{x}_{0}|}$ with $\nabla^2 G_{0}(\BS{x},\BS{x}_{0}) = - 4\pi \delta(\BS{x} - \BS{x}_{0})$, since $G_{k}(\BS{x},\BS{x}_{0}) \equiv G_{0}(\BS{x},\BS{x}_{0}) + \Delta G$ with 
\begin{align}
    \Delta G \equiv \frac{\exp{(\rmi k |\BS{x}-\BS{x}_{0}|)}-1}{|\BS{x}-\BS{x}_{0}|}
\end{align} 
which is regular when $\BS{x}\rightarrow \BS{x}_{0}$ \cite{YangJasa1997}. The same analysis and conclusion can be made for $\dff{G_{k}(\BS{x},\BS{x}_{0})}/\dff{n}$ (see Appendix~\ref{sec:app_notesDesing} for more details). Using this fact, we start with a known function $\psi(\BS{x})$ that satisfies the Laplace equation as, $\nabla^2 \psi(\BS{x}) = 0$, and the conventional boundary integral representation of it is 
\begin{align}
    c(\BS{x}_0) \psi(\BS{x}_{0}) + \int_{S} \psi(\BS{x}) \frac{\dff{G_{0}}}{\dff{n}} \text{ d} S(\BS{x}) = \int_{S}  \frac{\dff{\psi(\BS{x})}}{\dff{n}} G_{0} \text{ d} S(\BS{x}). \label{eq:cbimLapl}
\end{align}
Let us assume that $\psi(\BS{x})$ has the form 
\begin{align}
    \psi(\BS{x}) = g(\BS{x}) \phi(\BS{x}_0) + f(\BS{x}) \frac{\dff{\phi(\BS{x}_{0})}}{\dff{n}} \label{eq:nsbim_psi}
\end{align}
where $\phi(\boldsymbol x_0)$ and $\partial \phi/\partial n (\boldsymbol x_0)$ are constants in this context, and $g(\boldsymbol x)$ and $f(\boldsymbol x)$ satisfy the following conditions
\begin{subequations}\label{eq:nsbim_gfcond}
\begin{eqnarray} 
\nabla^2 g(\BS{x}) = 0, \quad \lim_{\boldsymbol x \to \boldsymbol x_0} g(\BS{x}) = 1, \quad \lim_{\boldsymbol x \to \boldsymbol x_0}\frac{\dff g(\BS{x})}{\dff{n}} = 0; \label{eq:nsbim_gcond}
\\
\nabla^2 f(\BS{x}) = 0, \quad \lim_{\boldsymbol x \to \boldsymbol x_0}f(\BS{x}) = 0, \quad \lim_{\boldsymbol x \to \boldsymbol x_0}\frac{\dff f(\BS{x})}{\dff{n}} = 1. \label{eq:nsbim_fcond}
\end{eqnarray}
\end{subequations}
Introducing Eq.~(\ref{eq:nsbim_psi}) into Eq.~(\ref{eq:cbimLapl}) and then subtracting the result from Eq.~(\ref{eq:BEM}):
\begin{align}
   & \int_{S} \left[\phi(\BS{x}) \frac{\dff{G_{k}}}{\dff{n}} - \phi(\BS{x}_0) g(\BS{x}) \frac{\dff{G_{0}}}{\dff{n}} + \phi(\BS{x}_0) \frac{\dff{g(\BS{x})}}{\dff{n}} G_{0} \right] \text{ d} S(\BS{x}) \nonumber \\
   & = \int_{S}  \left[ \frac{\dff{\phi(\BS{x})}}{\dff{n}} G_{k} - \frac{\dff{\phi(\BS{x}_0)}}{\dff{n}} \frac{\dff{f(\BS{x})}}{\dff{n}} G_{0} + \frac{\dff{\phi(\BS{x}_0)}}{\dff{n}} f(\BS{x}) \frac{\dff{G_{0}}}{\dff{n}} \right] \text{d} S(\BS{x}). \label{eq:nsbim}
\end{align}
Eq.~(\ref{eq:nsbim}) is the non-singular boundary integral equation for the Helmholtz equation, Eq.~(\ref{eq:Helm}), in which the integrands are all regular as $\BS{x}\rightarrow\BS{x}_0$ and the term $c(\boldsymbol x_0)\phi(\BS{x}_0)$ with the solid angle is eliminated. It is worth mentioning that Eq.~(\ref{eq:nsbim}) is valid for either $k=0$, $k$ being real, purely imaginary, or a random complex number. Note that each node $\boldsymbol x_0$ has its own $f$ and $g$ function. Further proof that Eq.~(\ref{eq:nsbim}) no longer contains singular terms can be found in Appendix~\ref{sec:app_notesDesing}.

There are many possible choices for functions $g(\BS{x})$ and $f(\BS{x})$ that can satisfy the conditions in Eq.~(\ref{eq:nsbim_gfcond}) which ensure Eq.~(\ref{eq:nsbim}) is free of singularities. For instance a constant and a linear function \cite{SunJCP2016}:
\begin{subequations} \label{eq:nsbim_gf}
%\begin{eqnarray} 
\begin{align}
g(\BS{x})  =  1 \qquad \qquad \qquad ,& \quad  \frac{\partial g (\boldsymbol x)}{\partial n}=0; \label{eq:nsbim_g}
\\
f(\BS{x})  =  \BS{n}(\BS{x}_0)\BS{\cdot}(\BS{x}-\BS{x}_0),& \quad \frac{\partial f (\boldsymbol x)}{\partial n}= \boldsymbol n (\boldsymbol x_0) \cdot \boldsymbol n (\boldsymbol x). \label{eq:nsbim_f}
\end{align}
%\end{eqnarray}
\end{subequations}
Note that for external problems, the integrals with the above choice of $\phi(\BS{x})$ with $g(\BS{x})$ and $f(\BS{x})$ in Eq.~(\ref{eq:nsbim_gf}) over the closed surface at infinity do not vanish. Nevertheless, the integral value can be found analytically as $4\pi \phi(\BS{x}_{0})$ and thus this value should be added to the left-hand side of Eq.~(\ref{eq:nsbim}) for external problems (see Appendix ~\ref{sec:app_inftyintegrals}).

It is worth emphasizing that the above framework enables us to use higher order surface elements, such as quadratic elements, together with the standard Gauss integration methods in the numerical implementation for all nodes (including the previously singular ones). In all the simulation results that follow, we have used the desingularized boundary element method as described above with quadratic elements, except for Fig.~\ref{fig:validation}(b) where we have used a desingularized version of the Burton-Miller boundary element method. The Sommerfeld radiation condition at infinity, i.e. only outgoing waves are allowed, is automatically satisfied (see also Appendix~\ref{sec:app_inftyintegrals}).

Assuming the surface is divided into surface elements with a total of $N$ nodes, then all potentials and their normal derivatives given in Eq.~(\ref{eq:nsbim}) are related by the following matrix system:
\begin{equation}
  {\cal{H}} \cdot \underline{\phi} = {\cal{G}} \cdot \underline{\partial \phi/\partial n}.
\end{equation}
The $N\times N$ matrices ${\cal{G}}$ and ${\cal{H}}$ are the numerical matrix equivalent of the boundary element integrals in Eq.~(\ref{eq:BEM}), and $\underline{\phi}$ and $\underline{\partial \phi/\partial n}$ are column vectors of length $N$ in which each component $\phi (i)$ and $\partial \phi(i)/\partial n$ belongs to $\phi$ and $\partial \phi/\partial n$, respectively, of the $i$'th node located at $\boldsymbol x_0 (i)$ with $i=1,...,N$. There are, by the way, other methods to solve the Helmholtz equation, such as with finite elements.

\subsection{Acoustic transducer}

\begin{figure}[t] 
\centering
\subfloat[]{\includegraphics[width=0.4\textwidth]{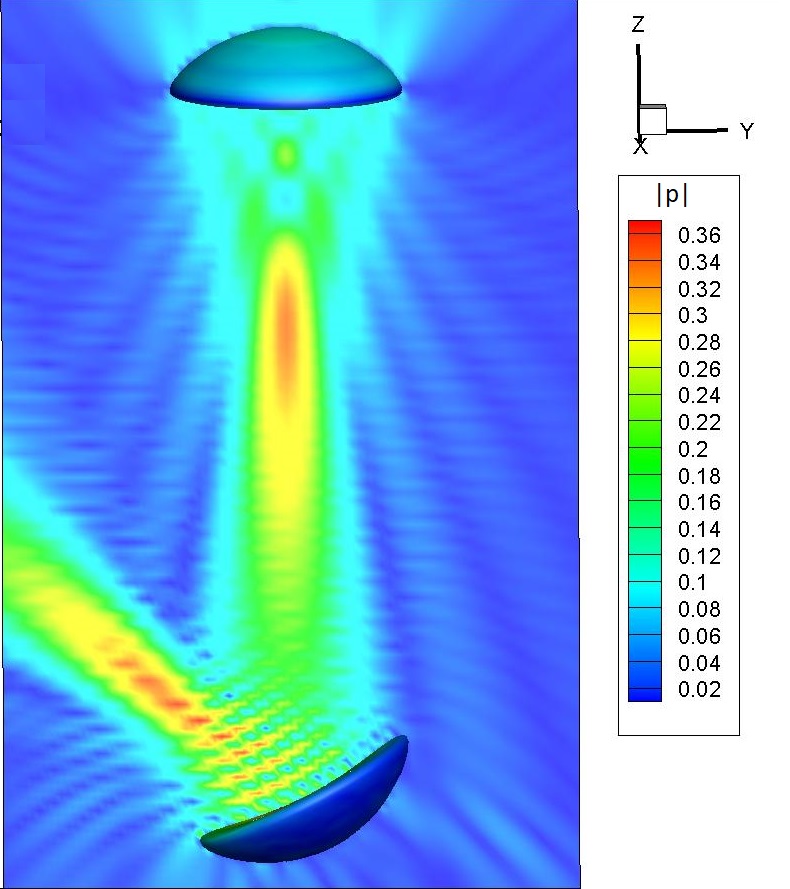}}\quad
\subfloat[]{\includegraphics[width=0.4\textwidth]{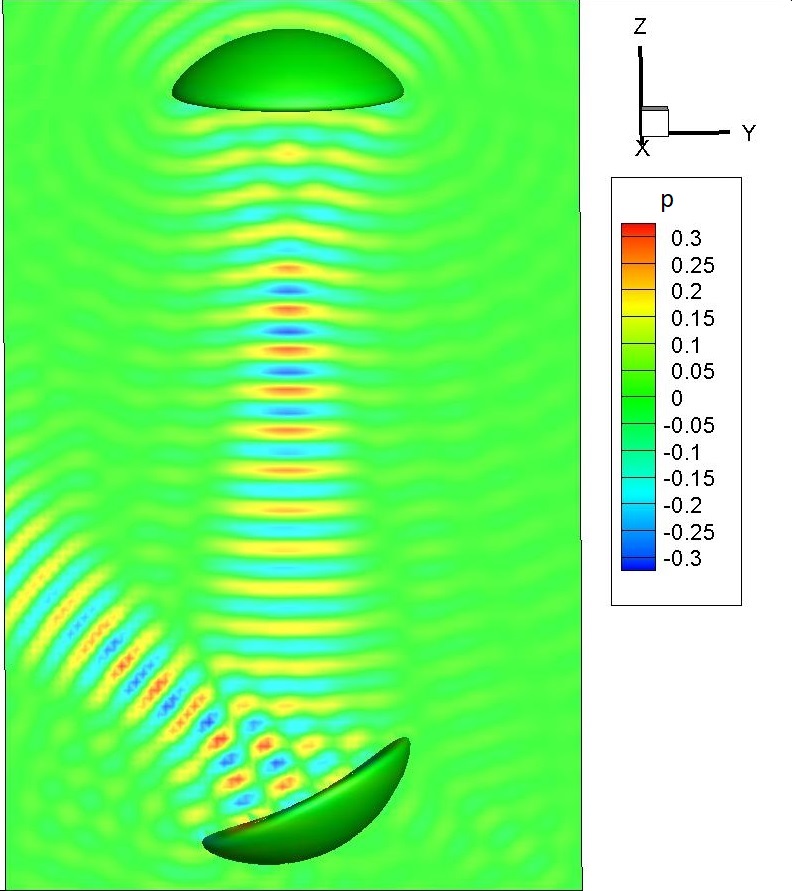}}
\caption{Sound transducer on top, receiving (stationary and rigid) rotated parabolic bowl below. (a) Absolute pressure amplitude and (b) instantaneous pressure profile. The graphs show the physical principles of generation and reflection of sound. Also note the focal areas of the two bowls. \label{fig:Transducer1}}
\end{figure}

To demonstrate some interesting wave phenomena in acoustics, in Fig.~\ref{fig:Transducer1}, we simulated an acoustic transducer (essentially a parabolic disk) oscillating up and down to generate sound waves which are then reflected onto a similar but stationary rigid disk that is rotated at an angle. Clearly, the focal points of the transducers can be observed. In Fig.~\ref{fig:Transducer1}(a), the absolute pressure is shown, while in Fig.~\ref{fig:Transducer1}(b), the instantaneous pressure profile is shown. When simulating this problem, only a surface mesh on the two disk boundaries is needed, and the pressure in the domain is obtained by postprocessing where the complex interference pattern is clearly visible.   

\subsection{The Helmholtz resonator} 

\begin{figure}[t] 
\centering
\subfloat[]{\includegraphics[width=0.2\textwidth]{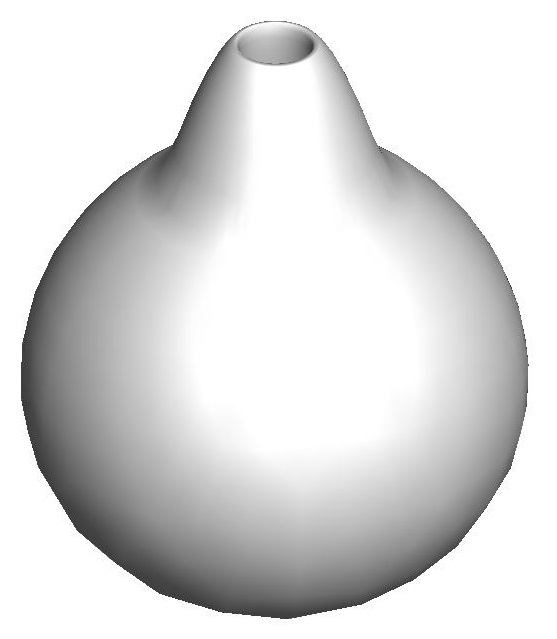}} \qquad
\subfloat[]{\includegraphics[width=0.2\textwidth]{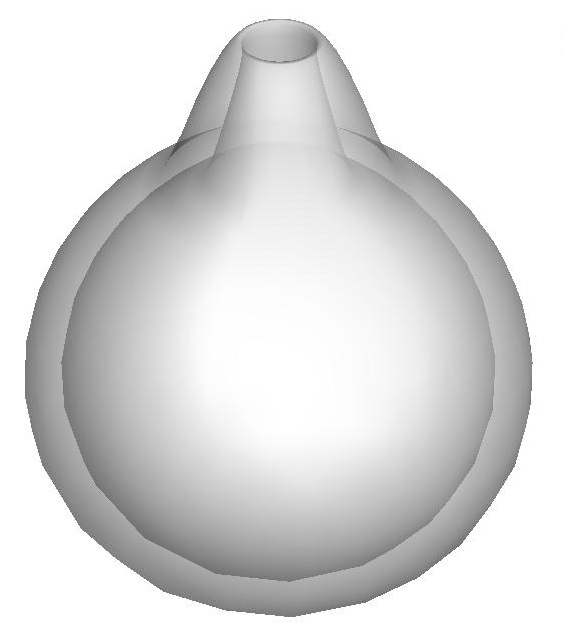}}
\caption{Helmholtz cavity as used in the simulations consisting of a central spherical part and a neck with an opening a) Full 3D view. b) Cross sectional view revealing the interior of the cavity. The theoretical resonance occurs at $ka=a \sqrt{A/VL}$, with $a$ the radius of the sphere, $A$ the opening area of the neck, $V$ the volume of the sphere and $L$ the length of the neck. For the current cavity this would be around $ka=0.225$.
\label{fig:Resonator1}}
\end{figure}

In this article, dedicated to Helmholtz, it seems appropriate to simulate an object that was named after him as well: a Helmholtz resonator. It consists of a thin spherical shell with an opening inside a neck at the top. The Helmholtz resonator is assumed to be acoustically hard: $\partial \phi^{\rmtot}/\partial n=0$ on its surface (thus $\partial \phi^{\rmsc}/\partial n=-\partial \phi^{\rminc}/\partial n$), where the superscripts $``\rmtot"$, $``\rminc"$ and $``\rmsc"$ indicate the total, incoming and scattered wave, respectively. In Fig.~\ref{fig:Resonator1} the Helmholtz cavity as used in our simulations is shown, both in full 3D view and in cross sectional view. The spherical part has an outer radius of $1.08a$, an inner radius of $0.92a$, a neck with length $L=0.67a$ and the neck tapers down from a value of $0.32a$ to $0.168a$ at the top. This will give us a theoretical resonance which should occur at $ka=a \sqrt{A/VL}=0.208$ with $A=(0.162a)^2 \pi$ the opening area of the neck, $V=4/3(0.92a)^3 \pi$ the inner volume of the sphere~\cite{RienstraHirschberg2016}. 

\begin{figure}[!t] 
\centering
\subfloat[]{\includegraphics[width=0.29\linewidth]{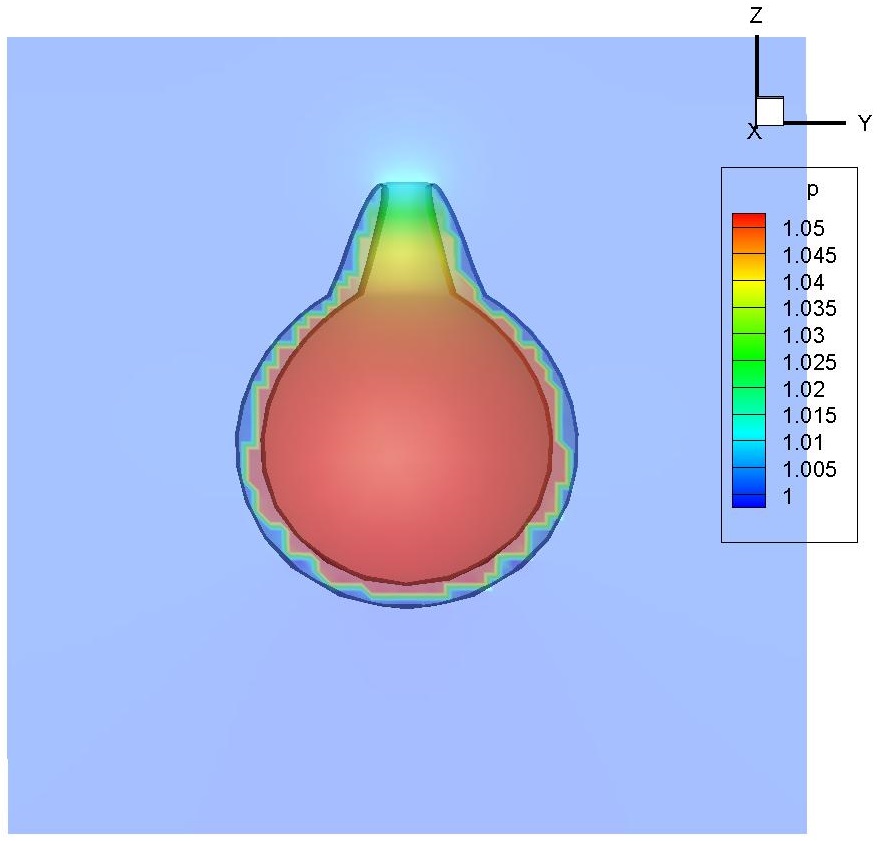}} \quad
\subfloat[]{\includegraphics[width=0.29\linewidth]{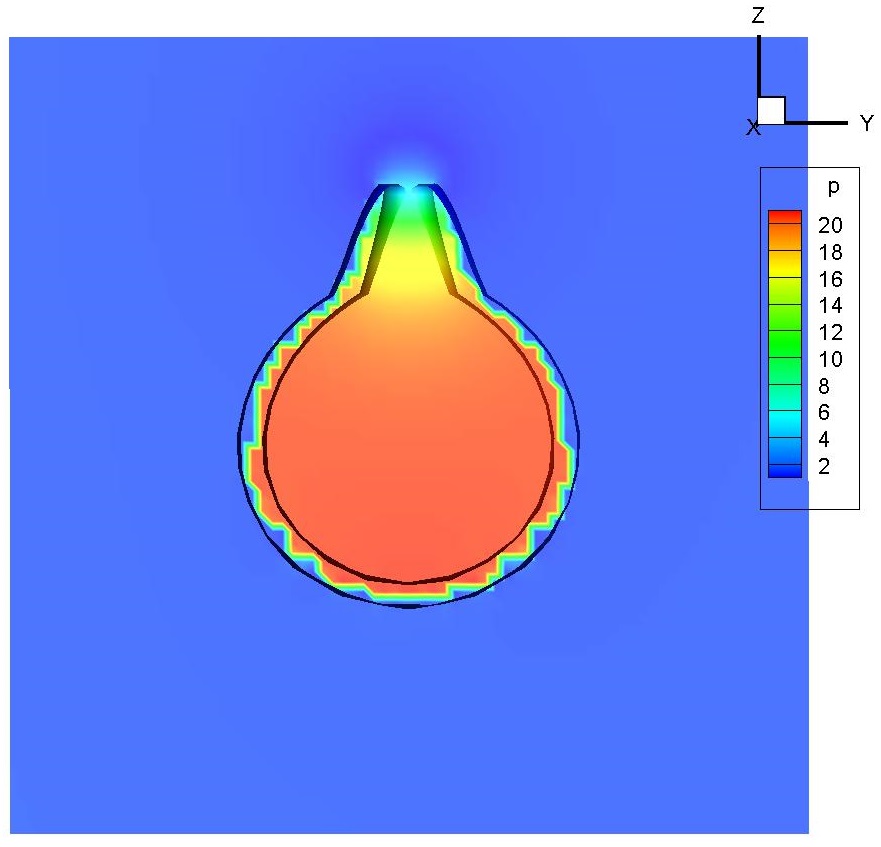}} \quad
\subfloat[]{\includegraphics[width=0.29\linewidth]{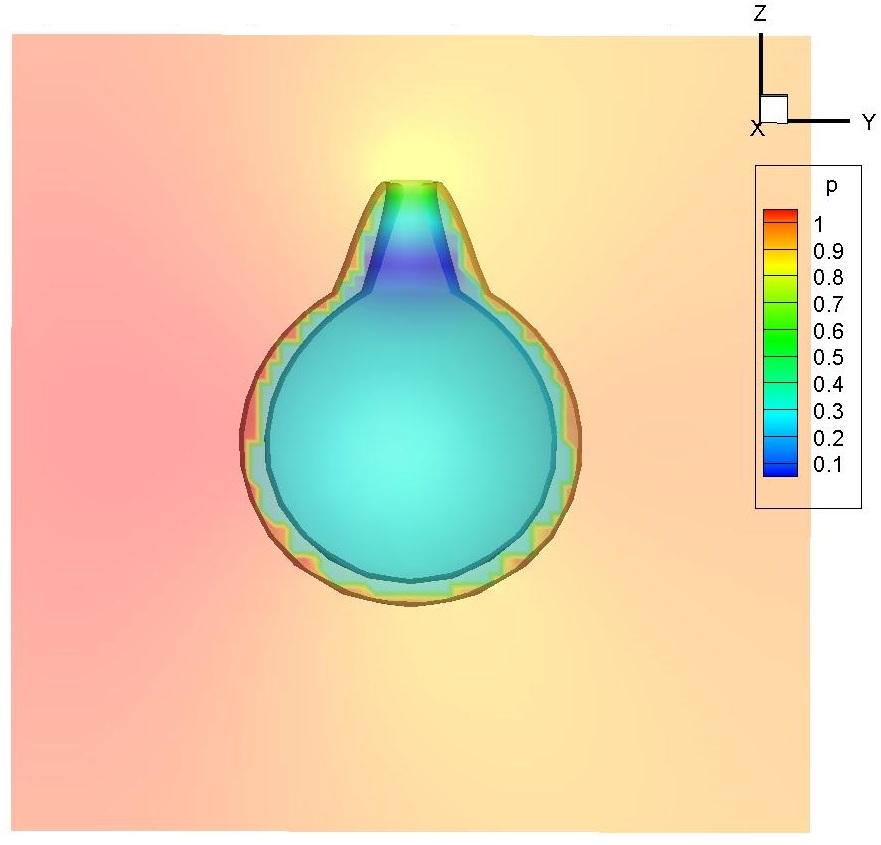}} \\
\subfloat[]{\centering\includegraphics[width=0.29\linewidth]{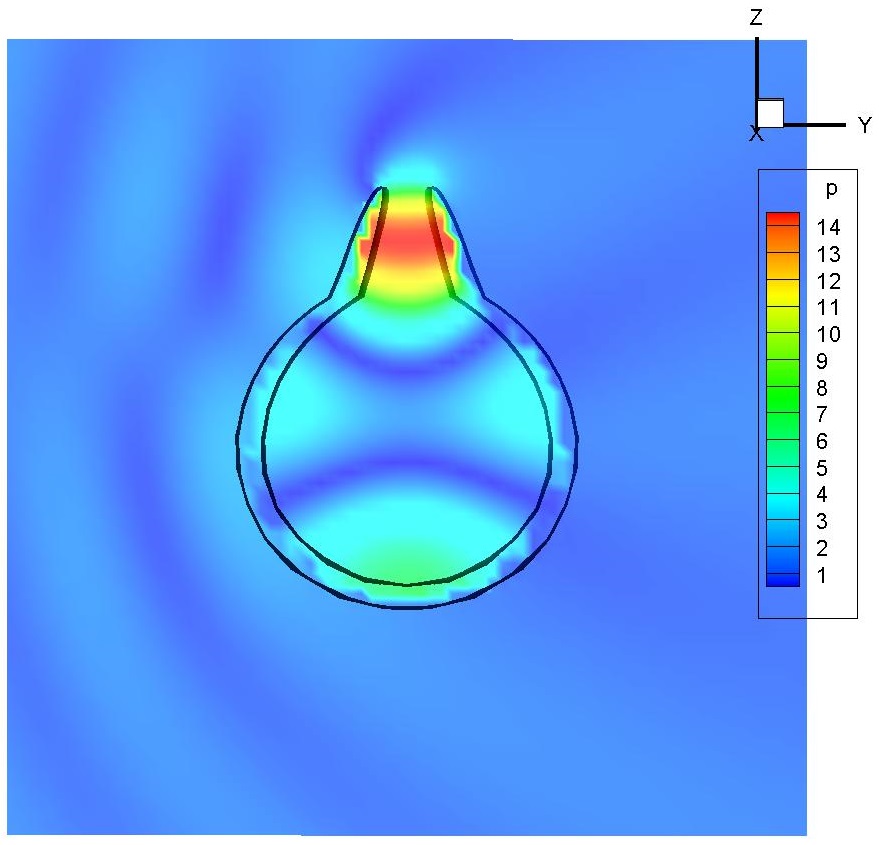}} \quad
\subfloat[]{\includegraphics[width=0.29\linewidth]{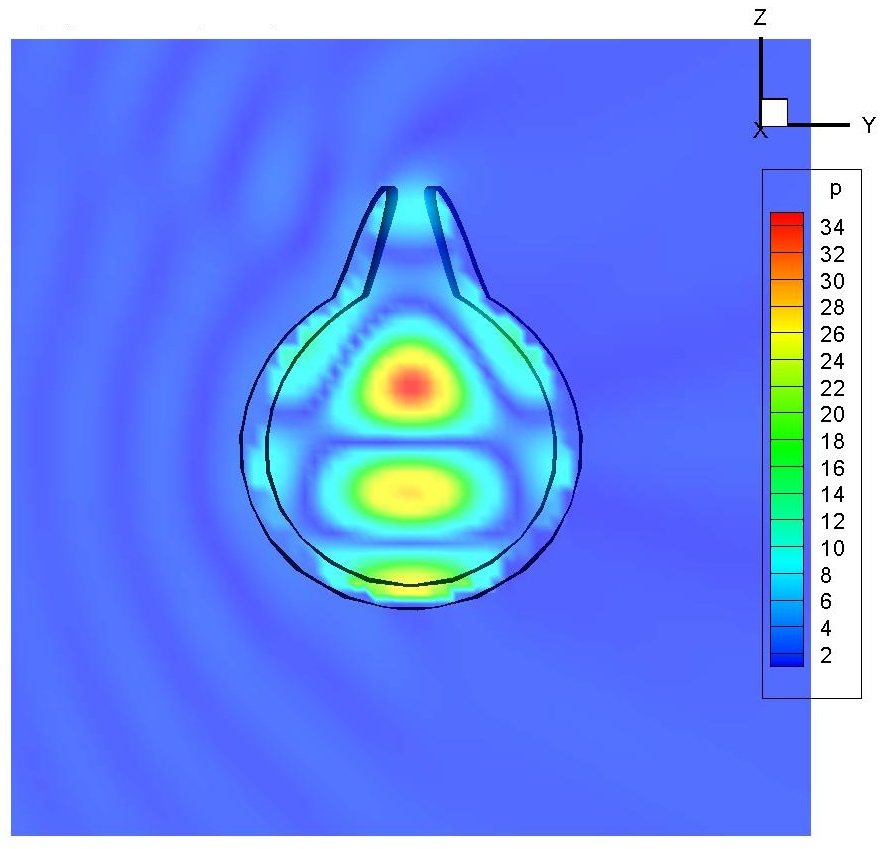}} \quad
\subfloat[]{\includegraphics[width=0.29\linewidth]{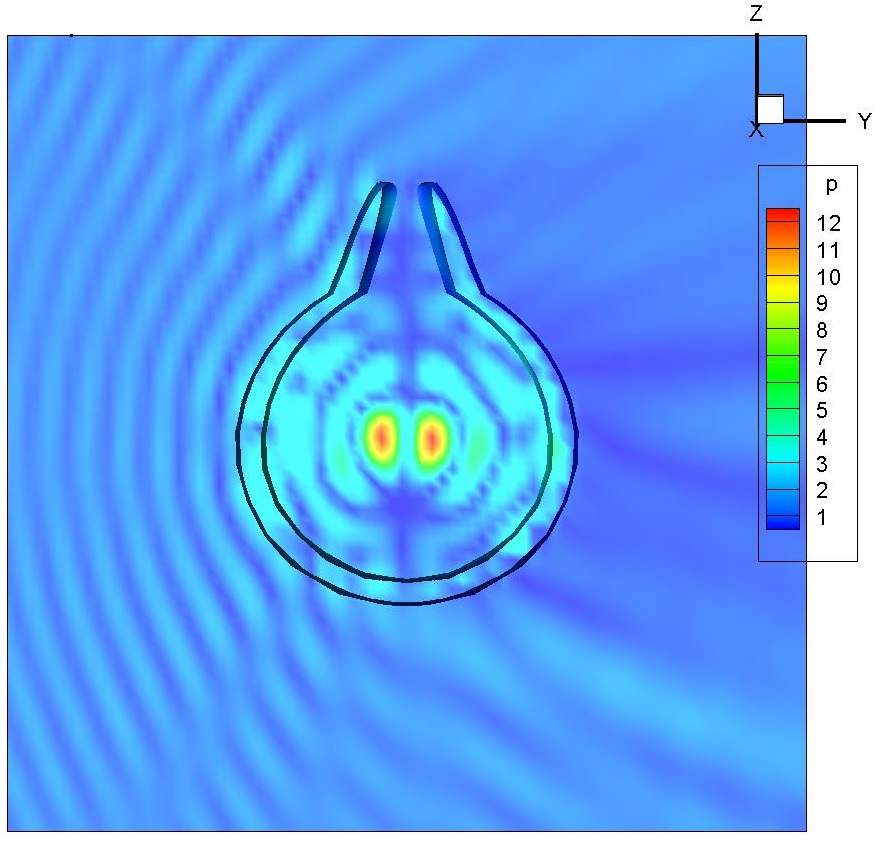}}
\caption{(a) $ka=0.05$, the pressure inside the cavity $p=1.05$ is slightly higher than the reference pressure outside $p=1.00$.  (b) $ka=0.23$ near the Helmholtz cavity resonance frequency.  (c) $ka=0.5$ The pressure $p=0.25$ is considerably lower than the reference pressure. (d) $ka=\pi$, near the first internal resonance frequency of a sphere. (e) $ka=2.05 \pi$ near the second resonance frequency of the sphere. (f) $ka=4.3\pi$, another internal resonance frequency.
\label{fig:Resonator2}}
\end{figure}

The incident acoustic wave travels from left to right in Fig.~\ref{fig:Resonator2}. We start with a very low $ka=0.05$ value in Fig.~\ref{fig:Resonator2}(a); the pressure inside the spherical part of the cavity is rather uniform, $p=1.05$ (slightly higher than the reference pressure of $p=1.00$). In Fig.~\ref{fig:Resonator2}(b) at $ka=0.23$, near the Helmholtz cavity resonance frequency, the pressure amplitude inside the resonator reaches a value of $p=20$ (i.e. 20 times the reference pressure). In Fig.~\ref{fig:Resonator2}(c) we increase $ka$ to a value of $ka=0.5$. The pressure is now $p=0.25$ and is considerably lower than the reference pressure. The first internal resonance frequency of the sphere is shown in Fig.~\ref{fig:Resonator2}(d) at $ka=\pi$. The maximum pressure is now situated in the neck and reaches a value of about $p=14$. A further resonance is shown in Fig.~\ref{fig:Resonator2}(e) for $ka=2.05 \pi$, near the second resonance frequency of the sphere. The maximum pressure amplitude reaches $p=34$ with multiple maximum values inside the sphere. Finally, in Fig.~\ref{fig:Resonator2}(f) yet another resonance is shown, now at $ka=4.3\pi$, with two horizontally placed pressure peaks of $p=12$. Also note the interference pattern to the left of the Helmholtz cavity for higher $ka$ values, which is caused by the interaction of the incoming and reflected wave on the external wall of the cavity.

\begin{figure}[t] 
\centering
{\includegraphics[width=0.6\textwidth]{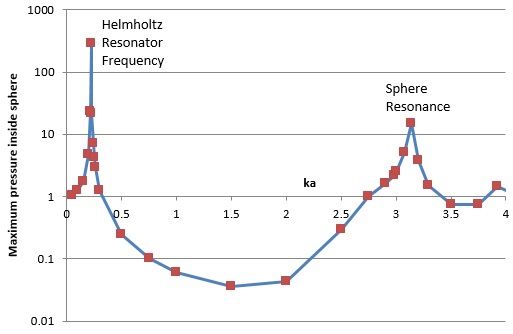}} 
\caption{The maximum pressure inside the Helmholtz cavity as a function of $ka$. A peak appears corresponding to the cavity resonance around $ka=0.225$ with a maximum of several 100's times the reference pressure. The peak at the internal resonance of the sphere is also visible around $ka=\pi$.  \label{fig:ResonatorExcel}}
\end{figure}

The maximum pressure in the Helmholtz cavity as a function of $ka$ is shown in Fig.~\ref{fig:ResonatorExcel}. Starting at frequency $ka=0$ the pressure inside the cavity is $1.0$ and thus equal to the reference pressure as it should be. It then reaches a very high peak near the Helmholtz cavity resonance around $ka=0.225$, this is close to the theoretically predicted value of $ka=0.208$. The difference can probably be attributed to the fact that the neck in our study is tapered and not straight. For large $ka$ values the pressure inside the cavity drops below the reference pressure and becomes rather low to reach a minimum around $ka=1.5$. A second peak can be observed near $ka=\pi$, which corresponds to the internal resonance frequency of a sphere. The maximum pressure does not necessarily occur in the main spherical part of the cavity, but can also occur in the neck part. Further peaks in the spectrum can be found associated with higher order internal resonance frequencies (not shown in Fig.~\ref{fig:ResonatorExcel}, but the pressure profiles are shown in Fig.~\ref{fig:Resonator2}(e) and (f) for some higher resonances near $ka=2\pi$ and $ka=4\pi$). The results clearly show that a Helmholtz cavity can be simulated, for both the low "Helmholtz resonance" as well as for the higher frequencies associated with the inner resonances of the spherical part.  

\section{Curl free vector Helmholtz equations} \label{sec:Curlfree} 

For a curl-free vector field, $\nabla \times \boldsymbol u = 0$ which also satisfies a Helmholtz equation $\nabla^2 \boldsymbol u + k^2 \boldsymbol u = \boldsymbol 0$, such as the case for the velocity field of a sound wave, we can introduce a potential $\boldsymbol u= \nabla \phi$. Then the framework of the previous section can be applied again. Thus curl free (or longitudinal) waves are relatively easy to describe using a single Helmholtz equation for the potential $\phi$. Such a simple approach is no longer available for divergence free vector fields as we will see in the next section.

\section{Divergence free vector Helmholtz equations; Electromagnetism}\label{sec:divfree} 
\subsection{Introduction to electromagnetic scattering}
Electromagnetic waves are typical examples of transverse waves. These waves satisfy the vector Helmholtz equation, but must simultaneously obey the zero divergence condition (a solenoidal vector field). Traditionally this has hindered the development of a simple and consistent framework just using Helmholtz equations alone and recourse had to be taken to Dyadic Green's functions or current based theories \cite{Kahnert2003}, for example the Poggio and Miller \cite{Poggio1973}, Chang and Harrington \cite{Chang1977}, and
Wu and Tsai \cite{Wu1977} (PMCHWT) theory, although alternative methods exist, such as the method of fundamental solutions \cite{Young2005,Leviatan1987}. In the frequency domain, if the medium is homogeneous and source free, the propagation of an electric field $\boldsymbol{E}$ obeys the following equations derived from the Maxwell's equations
~\cite{BohrenHuffman2008}
\begin{subequations}\label{eq:EM_Eqs}
\begin{eqnarray}
    \nabla^2 \BS{E} + k^2 \BS{E} &=& \BS{0}, \label{eq:Ecomp_Hohm} \\
    \nabla \cdot \BS{E} &=& 0. \label{eq:divE=0}
\end{eqnarray}
\end{subequations}
Eq.~(\ref{eq:divE=0}) is Gauss' law for electricity (assuming no volume charges). Eq.~(\ref{eq:Ecomp_Hohm}) is clearly a 3D version of the Helmholtz equation and is essentially composed of three scalar Helmholtz equations. In Eq.~(\ref{eq:Ecomp_Hohm}) (see pages 58-59 of \cite{BohrenHuffman2008}), the separate components of $\boldsymbol E$ only satisfy the scalar Helmholtz equation in a Cartesian coordinate system (i.e. $E_x$, $E_y$ and $E_z$). The difficulty encountered in electromagnetic wave theory is that the divergence free condition in Eq.~(\ref{eq:divE=0}) must be satisfied simultaneously with the vector wave equation~(\ref{eq:Ecomp_Hohm}). 

Obviously, there is the conceptual advantage of dealing directly with the physical quantity of interest, the electric field vector $\boldsymbol E$, instead of with surface currents or hypersingular dyadic Green's functions~\cite{Kahnert2003}. We will now introduce a recently developed field-only non-singular surface integral method to solve Eq.~(\ref{eq:EM_Eqs}) straightforwardly. In a typical simulation, we solve Eq.~(\ref{eq:EM_Eqs}) for the scattered electric field $\boldsymbol E^{\rmsc}$ of an incoming field $\boldsymbol E^{\rminc}$ in the external domain, since only the scattered field satisfies the Sommerfeld radiation condition at infinity. The total electric field would then be $\boldsymbol E^{\rmtot} = \boldsymbol E^{\rminc} + \boldsymbol E^{\rmsc}$, where $\boldsymbol E^{\rminc}$ is a plane wave in this work. If the object from which the scatter occurs is dielectric in nature, we also need to to solve for the transmitted electric field, $\boldsymbol E^{\rmtr}$ inside the object. 

\subsection{Boundary conditions for perfect electric conductors}  For the development of the numerical framework, it is convenient to express the field not only in global Cartesian components, but also in the local normal and tangential components on the scatterer surface (since the boundary conditions are often given in those components) as:
\begin{equation}
    \begin{aligned} \label{eq:E_conversion}
        \boldsymbol E = \boldsymbol e_x E_x + \boldsymbol e_y E_y + \boldsymbol e_z E_z 
        = \boldsymbol n E_n + \boldsymbol t_1 E_{t1} + \boldsymbol t_2 E_{t2}.
    \end{aligned}
\end{equation}
Here the normal vector, $\boldsymbol n$, is pointing out of the domain, and $\boldsymbol t_1$ and $\boldsymbol t_2$ are the two tangential vectors at the surface according to the convention $\boldsymbol n = \boldsymbol t_1 \times \boldsymbol t_2$ and $(\boldsymbol t_1 \cdot \boldsymbol t_2) = 0$. The unit vector in the $x$-direction is $\boldsymbol e_x$ (similar for $y$ and $z$). Thus for example the $x$-component of the electric field can be expressed in terms of the normal and tangential components as:
\begin{equation}\label{eq:ExDecomposition}
    E_x = (\boldsymbol n \cdot \boldsymbol e_x) E_n + (\boldsymbol t_1 \cdot \boldsymbol e_x) E_{t1} + (\boldsymbol t_2 \cdot \boldsymbol e_x) E_{t2}  = n_x E_n + t_{1x} E_{t1} + t_{2x} E_{t2}
\end{equation}
with $n_x = \boldsymbol n \cdot \boldsymbol e_x$, $t_{1x} = \boldsymbol t_1 \cdot \boldsymbol e_x$ and $t_{2x}= \boldsymbol t_2 \cdot \boldsymbol e_x$.

Perfect electric conductors are often  used as an idealization of a metallic object under the illumination of an electromagnetic wave. For a perfect electric conductor, the boundary condition is that the tangential component of the total electric field is zero. The total field is the sum of the incoming and scattered field thus:
\begin{align} \label{eq:Et1_BC}
E_{t1}^{\rmsc}=-E_{t1}^{\rminc} \quad ; \quad
E_{t2}^{\rmsc}=-E_{t2}^{\rminc} \quad \text{(Perfect Electric Conductor)}
\end{align}

\subsection{Divergence free condition implementation}

The key to solve Eq.~(\ref{eq:EM_Eqs}) is how to deal with the divergence free condition in Eq.~(\ref{eq:divE=0}). An elegant way to ensure that the divergence free condition of the electric field is satisfied in the domain, is to ensure it is satisfied on the surface of an object
. 
There are two things we need to investigate; firstly how do we ensure that the divergence free condition is satisfied on the surface and secondly, does that guarantee that the divergence free condition is satisfied everywhere else in the domain? These questions were answered in Sun et al.~\cite{SunJOSA_PEC2020} and will be summarized here. 

Let us start with the divergence of the electric field on the boundary of the domain:
\begin{equation} \label{eq:DivE_Surface}
    \nabla \cdot \boldsymbol E = \boldsymbol n \cdot \frac{\partial \boldsymbol E}{\partial n} - \kappa E_n + \frac{\partial E_{t1}}{\partial t_1} + \frac{\partial E_{t2}}{\partial t_2} =0 %\qquad \text{on surface}
\end{equation}
with $\kappa$ the curvature of the surface and $\partial/\partial n = \boldsymbol n \cdot \nabla$ the normal derivative, $\partial / \partial t_1 = \boldsymbol t_1 \cdot \nabla$, the tangential derivative in the tangential vector $\boldsymbol t_1$ direction and similar for $\partial / \partial t_2$ in the $\boldsymbol t_2$ direction. $E_n$, $E_{t1}$ and $E_{t2}$ are the normal and tangential components of the electric field on the scatterer surface, respectively. A convenient way to prove Eq.~(\ref{eq:DivE_Surface}) is to write the electric field into its normal and tangential components on the surface as $\boldsymbol E = E_n \boldsymbol n + E_{t1} \boldsymbol t_1 + E_{t2} \boldsymbol t_2$. Then the divergence can be written as $$\nabla \cdot \boldsymbol E = (\boldsymbol n \cdot \nabla) E_n  + E_n (\nabla \cdot \boldsymbol n) + (\boldsymbol t_1 \cdot \nabla) E_{t1}  + E_{t1} (\nabla \cdot \boldsymbol t_1) + (\boldsymbol t_2 \cdot \nabla) E_{t2}  + E_{t2} (\nabla \cdot \boldsymbol t_2).$$ With $-(\nabla \cdot \boldsymbol n) = \kappa$, the curvature, $(\nabla \cdot \boldsymbol t_1)=(\nabla \cdot \boldsymbol t_2)=0$, also $(\boldsymbol n \cdot \nabla) E_n = \frac{\partial }{\partial n} (\boldsymbol n \cdot \boldsymbol E)=\boldsymbol n \cdot \frac{\partial \boldsymbol E}{\partial n} + \frac{\partial n}{\partial n} \cdot \boldsymbol E$ and $\partial \boldsymbol n/\partial n = 0$, we can get back the desired equation, Eq.~(\ref{eq:DivE_Surface}).

Setting the divergence of the electric field on the boundary to zero is important since this will ensure that the divergence is also zero in the domain outside the object. This can be shown by realizing that if, for example $E_x$ is a solution of the Helmholtz equation, then $\partial E_x/ \partial x$ is as well. Then $\partial E_y/ \partial y$ and $\partial E_z/ \partial z$ are solutions as well. The sum of several solutions of the Helmholtz equation will also obey the Helmholtz equation, then $\nabla  \cdot \boldsymbol E$ will also obey the Helmholtz equation. Since $\nabla \cdot \boldsymbol E = 0$ on the surface, the normal derivative $ \partial (\nabla \cdot \boldsymbol E)/\partial n$ must also be zero (since these two quantities are related by the boundary element framework). Because the Helmholtz equation is elliptic in nature, the whole field must be divergence free. 

There are other ways to ensure the divergence free condition, for example using the vector identity $2(\nabla \cdot\BS{E}) \equiv \nabla^{2}(\BS{x}\cdot\BS{E}) + k^2 (\BS{x}\cdot\BS{E}) = 0$, with $\BS x$ the position vector. Thus one more Helmholtz equation for $(\BS x \cdot \BS E)$ will also guarantee $(\nabla \cdot \boldsymbol E) = 0$, see \cite{KlaseboerIEEETAP2017, SunPRB2017} for more details. 

The condition in Eq.~(\ref{eq:DivE_Surface}) can replace Eq.~(\ref{eq:divE=0}). On the surface of a perfect electric conductor, the total tangential fields $E_{t1}^{\rmtot}=0$ and $E_{t2}^{\rmtot}=0$, Eq.~(\ref{eq:DivE_Surface}) becomes:
\begin{equation} \label{eq:DivE_Surface2}
     \boldsymbol n \cdot \frac{\partial \boldsymbol E^{\rmtot}}{\partial n} = \kappa E_n^{\rmtot}.
\end{equation}
Thus for the scattered field we find:
\begin{equation} \label{eq:DivE_Surface3}
     \boldsymbol n \cdot \frac{\partial \boldsymbol E^{\rmsc}}{\partial n} =  \kappa E_n^{\rmsc} +\kappa E_n^{\rminc} - \boldsymbol n \cdot \frac{\partial \boldsymbol E^{\rminc}}{\partial n}.
\end{equation}
Similar to the decomposition of the electric field in Eq.~(\ref{eq:ExDecomposition}), a decomposition into normal and tangential components can be done for the normal derivatives as:
\begin{equation}
   \frac{\partial E_x}{\partial n} = n_x \left(\boldsymbol n \cdot \frac{\partial \boldsymbol E^{\rmsc}}{\partial n}\right) + t_{1x} \left(\boldsymbol t_1 \cdot \frac{\partial \boldsymbol E^{\rmsc}}{\partial n}\right) +t_{2x} \left(\boldsymbol t_2 \cdot \frac{\partial \boldsymbol E^{\rmsc}}{\partial n} \right)
\end{equation}
Thus the three matrix equations 
for $E_x^{\rmsc}$, $E_y^{\rmsc}$ and $E_z^{\rmsc}$ can be combined and expressed in terms of $E_n^{\rmsc}$, $ \BS{t}_{1} \cdot (\partial \BS{E}^{\rmsc}/\partial n)$ and $ \BS{t}_{2} \cdot (\partial \BS{E}^{\rmsc}/\partial n)$ for each node as a $3N \times 3N$ matrix system using Eq.~(\ref{eq:DivE_Surface3}):
  \begin{equation}
    \label{eq:linear_system_divE_Surface}
    \begin{aligned}
    &{\cal{H}} \cdot {E_x^{\rmsc}} = {\cal{G}} \cdot {\partial E_x^{\rmsc}/\partial n} \quad ; \quad 
    {\cal{H}} \cdot {E_y^{\rmsc}} = {\cal{G}} \cdot {\partial E_y^{\rmsc}/\partial n} \quad ; \quad 
    {\cal{H}} \cdot {E_z^{\rmsc}} = {\cal{G}} \cdot {\partial E_z^{\rmsc}/\partial n}\\
    &\qquad \qquad \qquad \Leftrightarrow
    \\
    &\begin{bmatrix}
      ({\cal{H}} - \kappa {\cal{G}})  n_x& - {\cal{G}} t_{1x} & -{\cal{G}} t_{2x}\\
       ({\cal{H}} - \kappa {\cal{G}}) n_y & -{\cal{G}} t_{1y} & -{\cal{G}} t_{2y} \\
      ({\cal{H}} - \kappa {\cal{G}}) n_z& -{\cal{G}} t_{1z}& -{\cal{G}} t_{2z}
    \end{bmatrix} \boldsymbol \cdot
    \begin{bmatrix}
      E_n^{\rmsc} \\ \BS{t}_{1} \cdot (\partial \BS{E}^{\rmsc}/\partial n) \\ \BS{t}_{2} \cdot (\partial \BS{E}^{\rmsc}/\partial n)
    \end{bmatrix}
    =
    \begin{bmatrix}
      {\cal{C}}_{x}\\ {\cal{C}}_{y} \\ {\cal{C}}_{z} 
    \end{bmatrix},
\end{aligned}
  \end{equation} 
  where ${\cal{C}}_x =  {\cal{H}} \left(t_{1x}  {E}^{\rminc}_{t_{1}}  +  t_{2x} {E}^{\rminc}_{t_{2}}  \right)  
  + n_x {\cal{G}}\left(\kappa {E}^{\rminc}_{n} - \BS{n} \cdot (\partial \BS{E}^{\rminc}/\partial n) \right)$, and similar for the $C_y$ and $C_z$ terms on the right hand side.
  It is also possible to work with the magnetic field instead of the electric field and get a similar $3N \times 3N$ matrix system. 

\begin{figure}[t] 
\centering\includegraphics[width=0.55\linewidth]{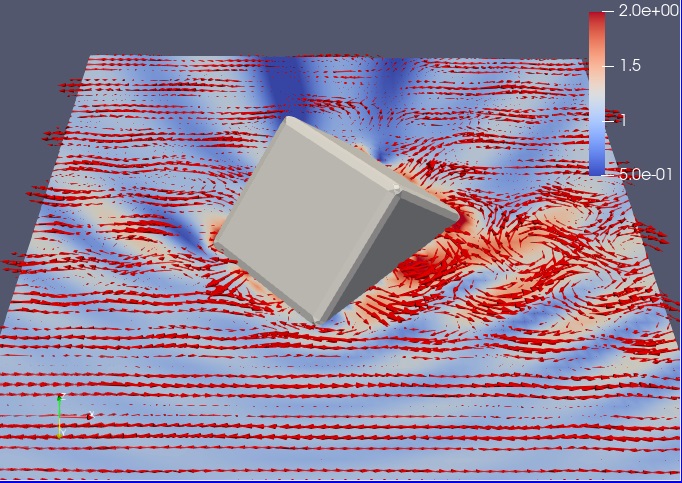}
\caption{Electromagnetic Scattering of a plane wave from a cube with rounded corners with length $2a$ and $ka=10$. The incoming plane wave travelling in the positive $z$-direction is polarized horizontally. The instantaneous electrical field vectors and amplitude are plotted on the plane $y=0$. The cube is rotated in such a way that one of its corners points towards the incoming wave. Note the strong reflection at the right side of the cube and the `shadow' region behind the cube. \label{fig:EMCube}}
\end{figure} 

The above formulation has been thoroughly tested against the Mie scattering solution for a sphere~\cite{SunJOSA_PEC2020}. Another example of scattering of an electromagnetic wave off a perfect conducting cube with rounded corners is shown in Fig.~\ref{fig:EMCube}. A cube with the edges having length $2a$ is rotated in such a way to have one vertex pointing in the direction of an incoming wave with $ka=10$. A complex pattern of interaction between the incoming and scattered waves can be seen in the plane at $y=0$. Since the right hand side exhibits a face, the perturbed electric field is rather chaotic there. The left hand side has an edge in the plane where we plot the electric field and the perturbation of the incoming field is much less there. Besides the instantaneous electric field (in red vectors), we have also plotted the absolute value of the electric field (which is time independent) as a scalar color plot.

\subsection{Dielectric formulation} 

Another type of scatterer that often appears is of dielectric nature. When a plane wave is scattered by a (closed) dielectric object, contrary to the perfect electric conductor case, we must also calculate the transmitted wave into the dielectric object. This means that we need to solve another set of Helmholtz equations for the dielectric with a different wave number $k_{\rmin}=\omega \sqrt{\epsilon_{\rmin} \mu_{\rmin}}$ as opposed to the external wave number $k_{\rmout}=\omega \sqrt{\epsilon_{\rmout} \mu_{\rmout}}$. To simplify the equations, we assume that for the permeabilities: $\mu_{\rmin}=\mu_{\rmout}$ (if these are different see Sun et al.~\cite{SunJOSA_Diel2020}), while the permitivities of both domains ($\epsilon_{\rmin}$ and $\epsilon_{\rmout}$) are different. 

We seek again a solution of the external scattered electric field $\boldsymbol E^{\rmsc}$ in terms of the incoming $\boldsymbol E^{\rminc}$ and transmitted field $\boldsymbol E^{\rmtr}$ into the dielectric object. The boundary conditions for such a system~\cite{Stratton1941} are (with no free surface charges or currents
):
\begin{subequations}
\begin{align} \label{eq:diel_BC_En}
    &\epsilon_{\rmin} E_n^{\rmtr} = \epsilon_{\rmout}(E_n^{\rminc} + E_n^{\rmsc}) \quad ; \quad \text{jump in normal electric field}\\
    \label{eq:diel_BC_Et1}
    &E_{t1}^{\rmtr} = E_{t1}^{\rminc} + E_{t1}^{\rmsc} \quad ; E_{t2}^{\rmtr} = E_{t2}^{\rminc} + E_{t2}^{\rmsc};\quad \text{Tang. continuity of electric field} \\
    \label{eq:diel_BC_Hn}
    &\mu_{\rmin} H_n^{\rmtr} = \mu_{\rmout}(H_n^{\rminc} + H_n^{\rmsc}) \quad ; \quad \text{jump in normal magnetic field}
    \\
    \label{eq:diel_BC_Ht1} 
    &H_{t1}^{\rmtr} = H_{t1}^{\rminc} + H_{t1}^{\rmsc} \quad ; H_{t2}^{\rmtr} = H_{t2}^{\rminc} + H_{t2}^{\rmsc} \quad \text{Tang. continuity of magnetic field} .
\end{align}
\end{subequations}
Since Eq.~(\ref{eq:diel_BC_En}) was derived from the divergence free electric field conditions on both sides of the boundary, we only have to ensure that $\nabla \cdot \boldsymbol E^{\rmsc} =0$ while $\nabla \cdot \boldsymbol E^{\rmtr}=0$ will then be automatically satisfied ($\nabla \cdot \boldsymbol E^\rminc=0$ already). The above boundary conditions are not all independent as Eq.~(\ref{eq:diel_BC_Et1}) guarantees that Eq.~(\ref{eq:diel_BC_Hn}) is satisfied.

The next task at hand is to convert the tangential continuity conditions on the magnetic field in Eq.~(\ref{eq:diel_BC_Ht1})  in terms of the electric field. Since $$\boldsymbol n \times \boldsymbol H = \boldsymbol n \times [\boldsymbol n H_n + \boldsymbol t_1 H_{t1} + \boldsymbol t_2 H_{t2}] =\boldsymbol n \times \boldsymbol t_1 H_{t1} + \boldsymbol n \times \boldsymbol t_2 H_{t2}  = \boldsymbol t_2 H_{t1} - \boldsymbol t_1 H_{t2}$$ this leads to $\boldsymbol t_2 \cdot [\boldsymbol n \times \boldsymbol H] = H_{t1}$ and $\boldsymbol t_1 \cdot [\boldsymbol n \times \boldsymbol H] = -H_{t2}$ (note the inversion of subscripts `1' and `2' here and the plus and minus signs). We can then get (by replacing $\boldsymbol H$ in $[\boldsymbol n \times \boldsymbol H]$, using Faraday's law of induction in the frequency domain: $\nabla \times \boldsymbol E = i \omega \mu \boldsymbol H$):
\begin{equation}
    \begin{aligned} \label{eq:diel_Ht1_Ht2}
        H_{t1} = \frac{1}{\rmi \omega \mu} \left[ \boldsymbol n \cdot \frac{\partial \boldsymbol E}{\partial t_2} - \boldsymbol t_2 \cdot \frac{\partial \boldsymbol E}{\partial n}\right] \quad ;\quad
        H_{t2} =- \frac{1}{\rmi \omega \mu} \left[ \boldsymbol n \cdot \frac{\partial \boldsymbol E}{\partial t_1} - \boldsymbol t_1 \cdot \frac{\partial \boldsymbol E}{\partial n}\right],
    \end{aligned}
\end{equation}
%
%In Einstein notation:
Since 
\begin{align}
\boldsymbol n \times \nabla \times \boldsymbol E = \epsilon_{ijk} n_j \epsilon_{klm} \partial E_m/\partial x_l = [\delta_{il} \delta_{jm} - \delta_{im} \delta_{jl}] n_j \partial E_m/ \partial x_l = n_j \partial E_j/\partial x_i - n_l \partial E_i/\partial x_l,
\end{align}
it will lead to 
\begin{subequations}
\begin{align}
    \boldsymbol t_1 \cdot (\boldsymbol n \times \nabla \times \boldsymbol E) &= \boldsymbol n \cdot \partial \boldsymbol E/ \partial t_1 - \boldsymbol t_1 \cdot \partial \boldsymbol E/ \partial n, \\
    \boldsymbol t_2 \cdot (\boldsymbol n \times \nabla \times \boldsymbol E) &= \boldsymbol n \cdot \partial \boldsymbol E/ \partial t_2 - \boldsymbol t_2 \cdot \partial \boldsymbol E/ \partial n,
\end{align}
\end{subequations}
and results in Eq.~(\ref{eq:diel_Ht1_Ht2}).

\begin{figure}[t] 
\centering
\subfloat[]{\includegraphics[width=0.41\textwidth]{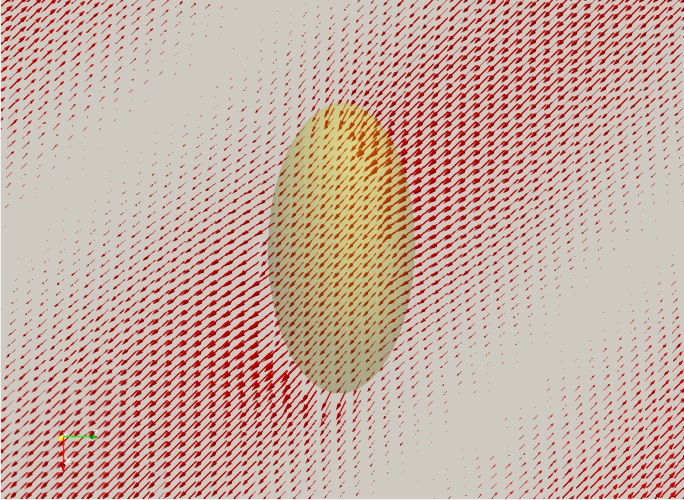}}\quad
\subfloat[]{\includegraphics[width=0.41\textwidth]{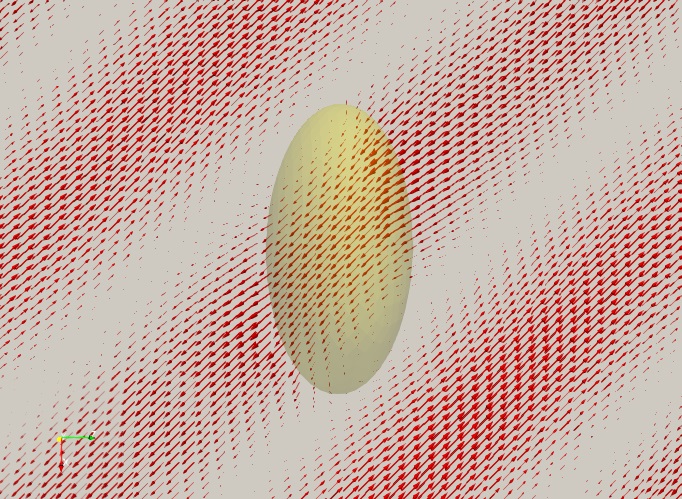}}
\caption{Instantaneous electric field vectors on the plane $y=0$ for a dielectric lens object with $a/b=2$ and the longest axis is $2a$ centered at $(x,y,z)=(0,0,0)$ under irradiation at an angle of 45 degrees (from top left to bottom right) (a) $k_{\rmin} a = 1.25$, $k_{\rmout}a=1.0$; (b) $k_{\rmin} a = 2.5$, $k_{\rmout}a=2.0$.   \label{fig:DielLens1}}
\end{figure}

The rest of the derivation for a dielectric scattering problem will result in a $6N \times 6N$ matrix system from which the quantities $E_n^{\rmsc}$, $E_{t1}^{\rmsc}$, $E_{t2}^{\rmsc}$, $\BS n \cdot \partial \BS E^{\rmsc}/ \partial n$, $\BS t_1 \cdot \partial \BS E^{\rmsc}/ \partial n$ and $\BS t_2 \cdot \partial \BS E^{\rmsc}/ \partial n$ are solved. From this the Cartesian components of the electric field and its normal derivatives can easily be reconstructed. The derivation is relatively straightforward and is given in Appendix~\ref{App:Dielectric}.

\begin{figure}[t] 
\centering
\subfloat[]{\includegraphics[width=0.41\textwidth]{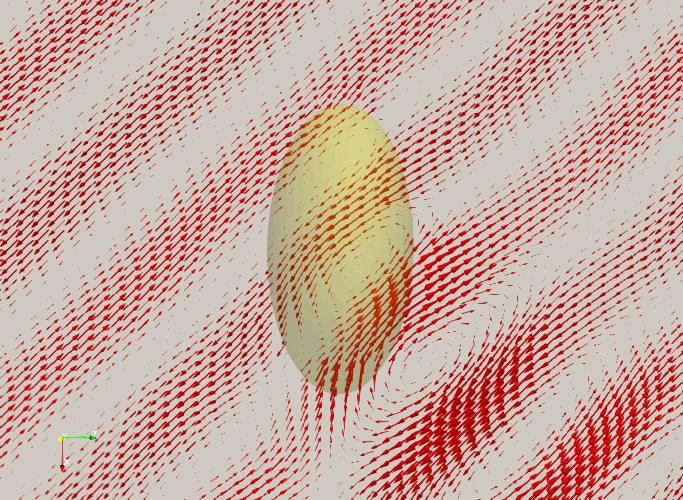}}\quad
\subfloat[]{\includegraphics[width=0.41\textwidth]{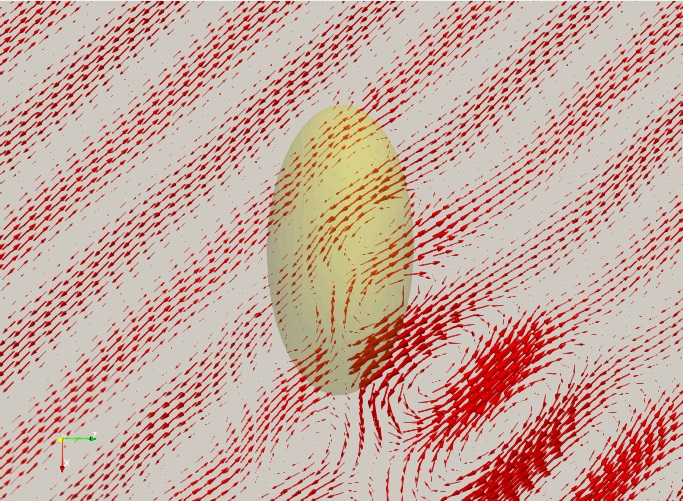}}
\caption{As Fig.~\ref{fig:DielLens1} but for higher $ka$ numbers:  (a) $k_{\rmin} a = 5$, $k_{\rmout}a=4$; (b) $k_{\rmin} a = 7.5$, $k_{\rmout}a=6$. Focusing effects start to appear behind the lens (in the lower right corner). 
\label{fig:DielLens2}}
\end{figure}

\begin{figure}[t] 
\centering
\subfloat[]{\includegraphics[width=0.41\textwidth]{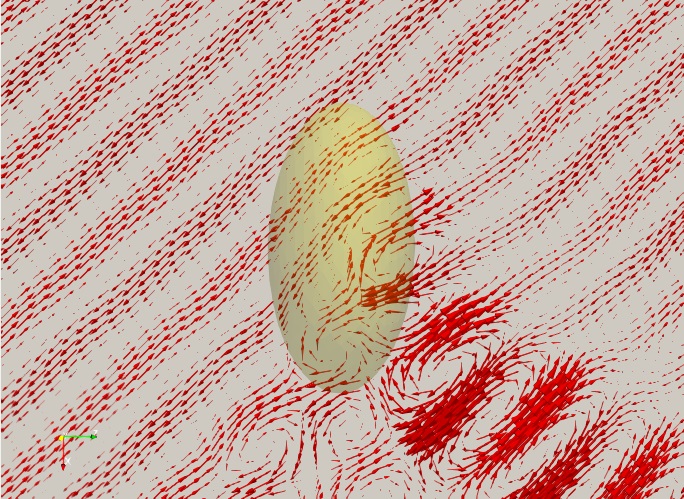}}\quad
\subfloat[]{\includegraphics[width=0.41\textwidth]{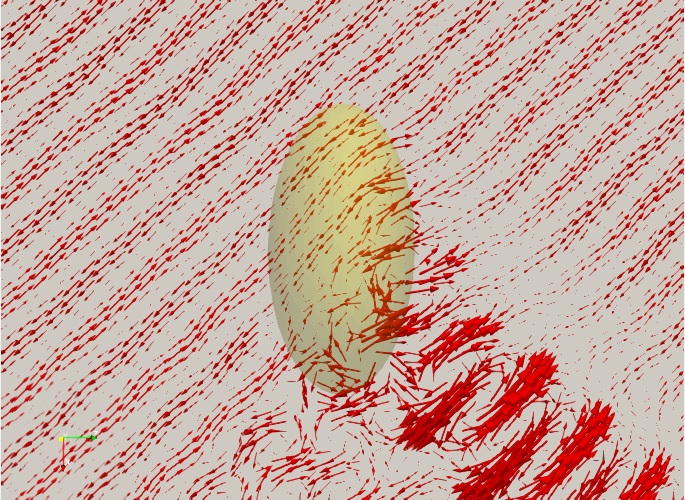}}
\caption{As Fig.~\ref{fig:DielLens1} but for even higher $ka$ numbers:  (a) $k_{\rmin} a = 10$, $k_{\rmout}a=8$ and (b) $k_{\rmin} a = 15$, $k_{\rmout}a=12$ (see also Fig.~\ref{fig:DielLens4}). The 'lens' gradually focuses the waves as the wavenumber gets higher.  \label{fig:DielLens3}}
\end{figure}
\begin{figure}[t] 
\centering
\subfloat[]{\includegraphics[width=0.45\textwidth]{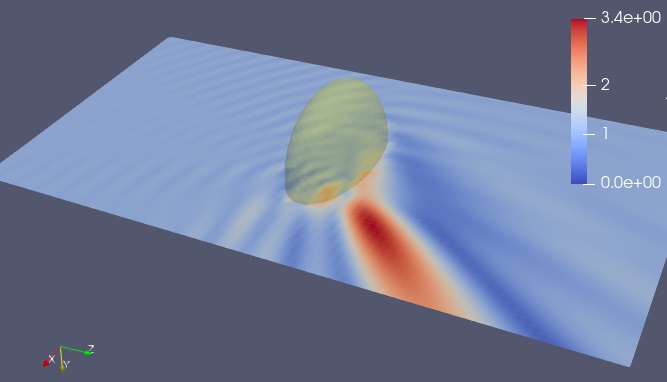}}\quad
\subfloat[]{\includegraphics[width=0.45\textwidth]{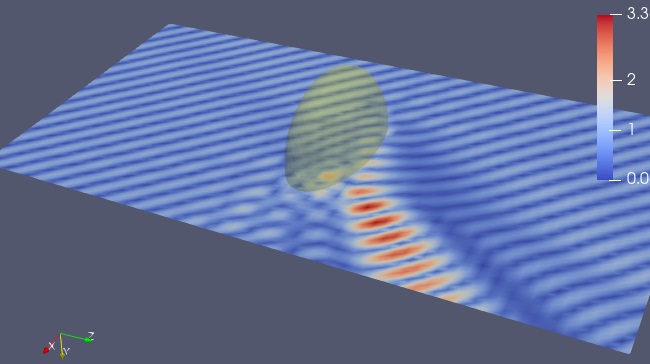}}\\
\subfloat[]{\includegraphics[width=0.45\textwidth]{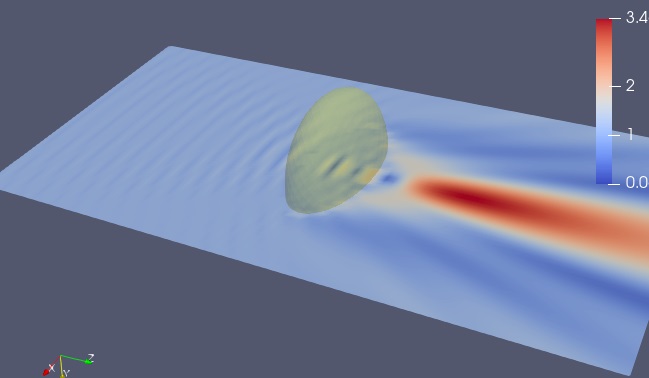}}\quad
\subfloat[]{\includegraphics[width=0.45\textwidth]{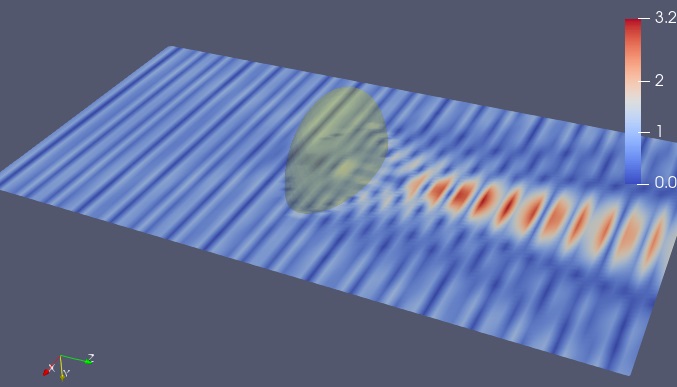}}
\caption{Electric field for a dielectric lens with $a/b=2$ under irradiation at an angle of 45 degrees with $k_{\rmin}a =15$ and $k_{\rmout}a =12$ (from top left to bottom right): (a) Absolute value of the electric field; (b) Absolute instantaneous value of the electric field. Under zero degrees irradiation, (c) Absolute value of the electric field; (d) Absolute instantaneous value of the electric field.  \label{fig:DielLens4}}
\end{figure}
As an illustrative example of a dielectric scattering problem, we consider an oblate spheroidal shaped object which functions as a lens, under illumination with a plane incoming wave at an angle with respect to the optical axis of the lens. The long axis of the lens has dimension $2a$ and the short axis is $a$. We first study the case when the incoming wave is at an angle of 45 degrees. For relatively low frequencies as shown in Fig~\ref{fig:DielLens1}, we can see that the electric field inside the lens is still more or less parallel. No focusing effect is observed for $k_{\rmin}a=1.25$ and $k_{\rmout}a=1.00$, neither for $k_{\rmin}a=2.5$ and $k_{\rmout}a=2.0$, which could be expected for the long wavelengths associated with these wavenumbers. We keep the material constants the same in all these examples (thus the ratio $k_{\rmin}/k_{\rmout}$ is fixed). Another example is shown next with $k_{\rmin}a=5$, $k_{\rmout}a=4$ in Fig.~\ref{fig:DielLens2}(a), and we can see some disturbance in the electric field in the lower right corner. This disturbance becomes bigger in Fig.~\ref{fig:DielLens2}(b) where we have used $k_{\rmin}a=7.5$ and $k_{\rmin}a=6$. In Fig.~\ref{fig:DielLens3}, we show two more examples for the parameters sets $k_{\rmin}a=10$, $k_{\rmout}a=8$ and $k_{\rmin}a=15$, $k_{\rmout}a=12$. We can see that the electric field vectors in the region to the right and below the lens are getting larger. In order to see more clearly what is going on here, we have plotted the absolute value of the complex valued electric field in Fig.~\ref{fig:DielLens4}(a). We can now clearly see the focal area of the lens. In Fig.~\ref{fig:DielLens4}(b) we have also plotted the length of the instantaneous electric field (thus the absolute value of the real part of the electric field vector), the incoming and focused waves can clearly be observed. As a reference case, we have also plotted these parameter for waves travelling at a zero angle towards the lens, Figs.~\ref{fig:DielLens4}(c) and (d).

\begin{figure}[t] 
\centering
\subfloat[]{\includegraphics[width=0.35\textwidth]{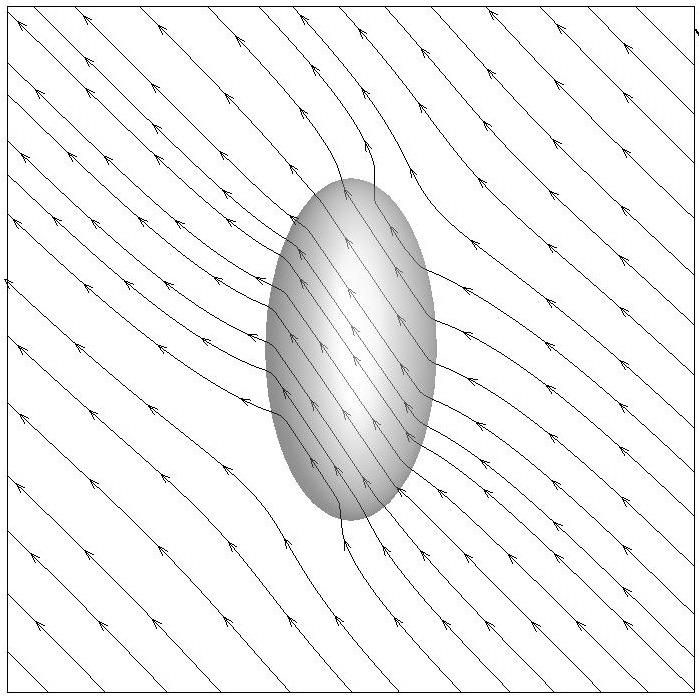}}\quad
\subfloat[]{\includegraphics[width=0.35\textwidth]{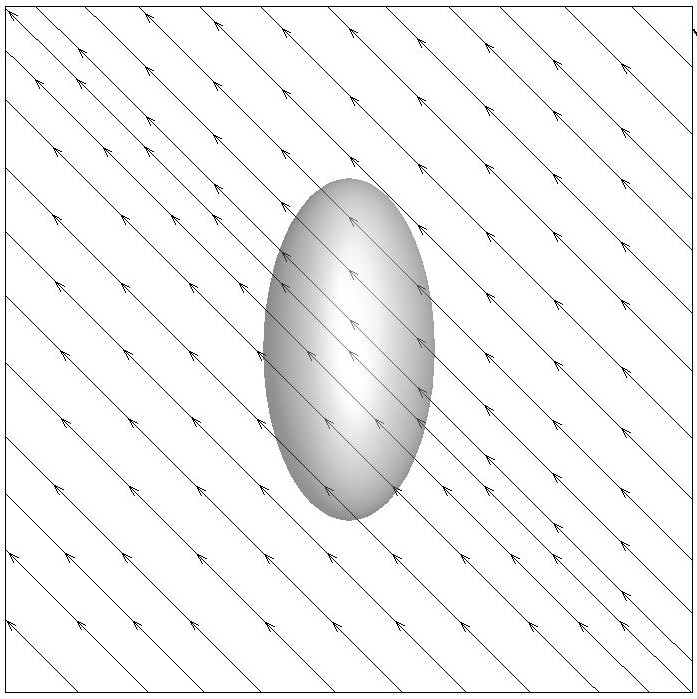}}
\caption{Electric field lines for a near zero frequency case ($k a \rightarrow 0$) for an ellipsoid dielectric object with longest length $2a$ under irradiation at an angle of 45 degrees. (a) With $k_{\rmout}a=0.001 $ and $k_{\rmin}a=0.002 $. Note that the internal electric field is parallel, but not parallel to the external field at infinity. (b) With $k_{\rmout}a=k_{\rmin}a=0.001$, the ellipsoid essentially becomes a transparent object. The examples clearly show that there is no 'zero frequency catastrophe' for this method. \label{fig:ZeroFreq1}}
\end{figure}

In Fig.~\ref{fig:ZeroFreq1}(a), the electric field lines around a dielectric object near the zero frequency limit are shown. Contrary to other numerical methods, such as surface current methods, our non-singular field only surface integral framework has no issues when $ka$ approaches zero (the long wavelength or electrostatic limit). We have chosen $k_{\rmin}a=0.002$ and $k_{\rmout}a=0.001$ for this particular example. Another example is shown in Fig.~\ref{fig:ZeroFreq1}(b), but now with $\epsilon_{\rmin} = \epsilon_{\rmout}$, thus the object essentially becomes invisible and the electric field lines are not being disturbed by the object. Here we have used $k_{\rmout}a=k_{\rmin}a=0.001$.

\section{Other physical systems with Helmholtz equations}
%\subsection{Dynamic Linear Elasticity}\label{sec:DLE}
One physical phenomenon of special scientific interest concerns acoustic waves in solid materials since transverse and longitudinal waves, each with a different wave numbers, can occur simultaneously. The transverse waves satisfy a zero divergence, while for the longitudinal waves, the curl is zero. These waves travel at different speeds (the longitudinal waves always travel faster than the transverse waves). Assume a linear isotropic homogeneous elastic material, with zero body force. The material has density $\rho$ and the Lam\'e constants $\lambda$ and $\mu$ ($\mu$ is the shear modulus) as material elastic properties. 

In the time domain waves occurring in the material can be described with the following elastodynamic equation
$\rho \frac{\partial^2 \boldsymbol{u'}}{\partial t^2} = \nabla \cdot \boldsymbol{\sigma}$ where $
    \boldsymbol{\sigma} = \lambda  (\nabla \cdot \boldsymbol{u'}) \textbf{I} + \mu [\nabla \boldsymbol{u'} + (\nabla \boldsymbol{u'})^T ]$,
with $\boldsymbol u'$ the displacement vector, and $\boldsymbol \sigma$ the stress tensor. The superscript $T$ indicates the transpose of the tensor and $\boldsymbol I$ is the unit tensor. This equation essentially expresses the force balance on an infinitesimal volume element; with inertial forces on the left and elastic forces on the right. In the frequency domain this leads to the Navier equation (assuming again a harmonic time dependency $\boldsymbol u'=\boldsymbol u e^{-\rmi \omega t}$)~\cite{Eringen1975,Beskos1987}:
$ c_L^2 \nabla \nabla \cdot \boldsymbol u
   - c_T^2 \nabla \times \nabla \times \boldsymbol u + \omega^2 \boldsymbol u = \boldsymbol 0$
where $c_T^2 = \mu/\rho$ and 
$c_L^2 = (\lambda +2 \mu)/\rho$.  
$c_T$ and $c_L$ are the transverse and longitudinal wave speeds, also often referred to as shear wave velocity and dilatational wave velocity. Alternatively, using the vector identity $\nabla \times \nabla \times \boldsymbol u = \nabla \nabla \cdot \boldsymbol u - \nabla^2 \boldsymbol u$, the Navier equation can be written as: 
\begin{equation} \label{eq:Navier2}
  \left( \frac{k_T^2}{k_L^2}-1 \right)  \nabla \nabla \cdot \boldsymbol u
   + \nabla^2  \boldsymbol u + k_T^2 \boldsymbol u = \boldsymbol 0
\end{equation}
with $k_L=\omega/c_L$ and $k_T=\omega/c_T$ the longitudinal and transverse wave numbers, respectively. We can now describe the waves occurring in a linear elastic medium with the help of Helmholtz equations alone, using the Helmholtz decomposition. The displacement field is decomposed into longitudinal and transverse components $\boldsymbol{u} = \boldsymbol{u}_T + \boldsymbol{u}_L$ where $\nabla \cdot \boldsymbol{u}_T = 0$ (transverse) and $\nabla \times \boldsymbol{u}_L = \boldsymbol{0}$ (longitudinal)~\cite{KlaseboerJE2018}, then Eq.~(\ref{eq:Navier2}) becomes:
\begin{equation}
\begin{aligned} 
    \nabla^2 \boldsymbol{u}_T + k_T^2 \boldsymbol{u}_T = \boldsymbol{0} \quad &; \quad
    \nabla \cdot \boldsymbol u_T = 0 \\
    \nabla^2 \Phi + k_L^2 \Phi = 0 \quad &.  \label{eq:x_cross_uL}
\end{aligned}
\end{equation}
Here we have chosen a potential function $\boldsymbol u_L=\nabla \Phi$ to automatically satisfy the curl-free condition for $\boldsymbol u_L$. The above framework was used recently to find an analytical solution for a vibrating rigid sphere (i.e. moving periodically up and down in an infinite elastic material)~\cite{KlaseboerJE2018} and for the same sphere but surrounded by an additional elastic shell~\cite{KlaseboerElasticShell2022}. 

The Helmholtz equation is classically used to describe acoustic waves, and the wave number $k$ is often a real number. However, $k$ can be complex, and Eq.~(\ref{eq:Helm}) then represents an acoustic wave in a medium with absorption. As an example we can mention the description of acoustic boundary layers around a sphere in a viscous liquid \cite{KlaseboerPhysFl2020}, which was based on the Navier equation, Eq.(\ref{eq:Navier2}), but with complex $k$'s.
There are other types of variations of the Helmholtz equation as well. For example, when $k = i \kappa_{\text{DH}}$ is imaginary (thus $k^2$ is negative), Eq.~(\ref{eq:Helm}) becomes the Debye-H\"uckel model for the molecular electrostatics potential $\Psi$, satisfying $
   \nabla^2 \Psi - \kappa_{\text{DH}}^2 \Psi = 0
$, in colloidal systems in which $\kappa_{\text{DH}}$ is the inverse of the Debye length~\cite{SunJCP2016}. When $k^2$ is purely imaginary, Eq.~(\ref{eq:Helm}) describes diffusion or heat transfer and no longer wave phenomena. Nevertheless, the same numerical boundary element framework (see Sec.~\ref{sec:BEMChapter}) can still be applied.

\section{Discussion} \label{sec:Discussion}
All the above classic problems can, in principle, be solved just based on one or a set of Helmholtz equations. The following points are worth noting concerning the non-singular boundary element framework of Section~\ref{sec:BEMChapter}.

\begin{itemize}
\item Since we use the same numerical framework for acoustics and electromagnetics, we can use the same codes with quadratic elements for both. For the electromagnetic problem we do not use RWG~\cite{RWS1982} elements. Thus, anyone with a Helmholtz boundary element solver for sound waves, can in principle extend this to electromagnetic scattering simulations. Moreover, our method is stable for electromagnetic problems in the long wavelength limit. 

\item When solving vector wave problems, such as in electromagnetics, iterative solvers can also be used only requiring one $N \times N$ matrix to be solved, although careful attention need to be paid to the convergence behavior.
\item In Fig.~\ref{fig:validation}(a), the convergence behavior of the solution of the Helmholtz equation is given for scattering on a hard sphere for $ka=1.01\pi$. The solution at a sample point at $r=1.2a$ (see the inset of the figure) is compared to the theoretical solution and the error is plotted as a function of node number. The relative error is about $10^{-3}$ for 1000 nodes and goes down to $10^{-5}$ for about 8000 nodes. The required CPU time is also indicated in the same graph. 
\item Concerning problems in the time domain: a time dependent problem can always be decomposed into its Fourier components. Then those components can be treated in the frequency domain and the solution can be reassembled in the time domain. This was actually implemented for sound waves and for electromagnetic pulses in Klaseboer et al.~\cite{KlaseboerJASA2017,KlaseboerAO2017}.

\end{itemize}

When using a boundary element method to numerically solve the Helmholtz equation for external problems, one challenge is that non-physical solutions will show at certain discrete frequencies~\cite{MarburgWu2008, Misawa2017}. This numerical issue is known as the occurrence of spurious solutions and occurs at  fictitious frequencies. The two most popular methods to deal with fictitious frequencies are the CHIEF method proposed by Schenck~\cite{Schenck1968} and the Burton-Miller method~\cite{BurtonMiller1971}. The CHIEF method uses additional internal points creating an over-determined system. Therefore an additional computational technique, such as the least square method must be employed. Furthermore, the successful removal of spurious solutions at fictitious frequencies cannot be fully guaranteed by this method. The Burton-Miller method claims that the spurious solutions at the fictitious frequencies disappear. However, the drawback is that it contains an integral with a hyper-singularity which involves integration in parts and a principle value integration. 
Despite the fact that the Burton-Miller formulation results in hyper singular integrals, it is still possible to fully desingularize all the integrals involved. Since we believe this fully non-singular Burton-Miller formulation has not been shown before, the derivation is given in detail in Appendix~\ref{App:nsbim_BM}. In Fig.~\ref{fig:validation}(b) the scattering from a rigid sphere is shown where the desingulared Burton-Miller framework of Eq.~(\ref{eq:nsbimBM_a}) is used for $ka=0$ to $30$. As expected, no fictitious frequencies appear here.

\section{Conclusions} \label{sec:Conclusion}

In this article, in memory of the famous German scientist Hermann von Helmholtz for his 200-year birthday, we revisited the wide engineering applications of the elliptic partial differential equation named after him, the Helmholtz equation, and introduced a robust, efficient and simple non-singular boundary integral (element) method to solve the Helmholtz equation. We have shown that the Helmholtz equation, classically used to describe scalar quantities, such as the potential or pressure in sound wave simulations in the frequency domain, can also be used for other physical systems. We have given an example for the case of electromagnetic waves, which is essentially a set of three coupled Helmholtz equations, one for each Cartesian component of the electric field. The main difficulty that is faced historically, is how to enforce the divergence free condition of the electric field. We have shown that this condition can be satisfied relatively easily. Some examples are shown for perfect electric conductors and dielectric objects. For the scalar Helmholtz equation we have chosen an example with transducers and reflectors and a Helmholtz cavity.

A major advantage in solving actual problems is the development of an entirely desingularized boundary element method for the Helmholtz equation without compromising on accuracy or efficacy. The main idea behind this being that if the physical problem does not have any singularities, then the mathematical model also should not have any singular behavior. The singular behavior in the classical boundary element method originates from the Green's function. It turns out that if a carefully chosen analytical solution, with the same singular behavior, is subtracted from the original equation, a totally non singular framework is created. All elements, including the previously singular ones, can be integrated with standard Gauss integration. It also enables us to implement higher order quadratic elements with quadratic shape functions with ease. 

\appendix

\section{Notes on non-singular boundary element method}\label{sec:app_notesDesing} 

A closer inspection of the desingularized equation (\ref{eq:nsbim}) will be performed here. The first and second terms in Eq.~(\ref{eq:nsbim}), i.e. $\phi(\boldsymbol x) \frac{\partial G_k}{\partial n} - \phi(\boldsymbol x_0) g(\boldsymbol x) \frac{\partial G_0}{\partial n}$ contain no singularities as can be shown by taking a Taylor expansion at $\boldsymbol x_0$ as: 
\begin{equation}
    \begin{aligned}\label{eq:AppBTaylor}
    &\phi(\boldsymbol x) \frac{\partial G_k}{\partial n} - \phi(\boldsymbol x_0) g(\boldsymbol x) \frac{\partial G_0}{\partial n}\\
    \approx& \phi(\boldsymbol x_0) \frac{\partial G_k}{\partial n} +\nabla \phi \cdot (\boldsymbol x - \boldsymbol x_0) \frac{\partial G_k}{\partial n} - \phi(\boldsymbol x_0) \frac{\partial G_0}{\partial n}- \phi(\boldsymbol x_0) \nabla g \cdot(\boldsymbol x-\boldsymbol x_0) \frac{\partial G_0}{\partial n}\\
    =&\phi(\boldsymbol x_0)\left[ \frac{\partial G_k}{\partial n}-\frac{\partial G_0}{\partial n} \right] +\nabla \phi \cdot (\boldsymbol x - \boldsymbol x_0) \frac{\partial G_k}{\partial n} - - \phi(\boldsymbol x_0) \nabla g \cdot(\boldsymbol x-\boldsymbol x_0) \frac{\partial G_0}{\partial n}.
    \end{aligned}
\end{equation}
The expressions for $\partial G_g/\partial n$ and $\partial G_0/ \partial n$ are (with $r=|\boldsymbol x - \boldsymbol x_0|$):
\begin{equation} \label{eq:AppBdGdn}
    \begin{aligned} 
    \frac{\partial G_k}{\partial n} &=[\rmi kr-1]\frac{e^{\rmi kr}}{r^3}(\boldsymbol x - \boldsymbol x_0) \cdot \boldsymbol n(\boldsymbol x) 
    \quad ; \quad
    \frac{\partial G_0}{\partial n} &=-\frac{1}{r^3}(\boldsymbol x - \boldsymbol x_0) \cdot \boldsymbol n(\boldsymbol x) 
    \end{aligned}
\end{equation}
The terms in between brackets in Eq.~(\ref{eq:AppBTaylor}) can thus be written with Eq.~(\ref{eq:AppBdGdn}) as:
\begin{equation} \label{eq:AppBbracket}
\begin{aligned}
\left[ \frac{\partial G_k}{\partial n}-\frac{\partial G_0}{\partial n} \right]
=\left([\rmi kr-1]e^{\rmi kr}+1\right)\frac{1}{r^3}(\boldsymbol x - \boldsymbol x_0) \cdot \boldsymbol n(\boldsymbol x).
\end{aligned}
\end{equation}
Here $[\rmi kr-1]e^{\rmi kr}+1= 2 \rmi kr + o(k^2r^2)$ and the term $(\boldsymbol x - \boldsymbol x_0) \cdot \boldsymbol n (\boldsymbol x)$ behaves as $1/r^2$ since in the limit of $\boldsymbol x$ going to $\boldsymbol x_0$, $(\boldsymbol x - \boldsymbol x_0)$ and $\boldsymbol n$ are perpendicular. Thus Eq.~(\ref{eq:AppBbracket}) does not contain any singular term.  The term with $\nabla \phi \cdot (\boldsymbol x - \boldsymbol x_0)$ in Eq.~(\ref{eq:AppBTaylor}) goes as $r$ and thus cancels out the $1/r$ behavior of $\partial G_k/ \partial n$. Similar for $\nabla g \cdot (\boldsymbol x - \boldsymbol x_0)$. 

For the third term in Eq.~(\ref{eq:nsbim}) with $\phi(\boldsymbol x_0) \frac{\partial g(\boldsymbol x)}{\partial n} G_0$, as long as $\partial g/\partial n$ approaches zero in a linear manner, it will cancel out the $1/r$ singularity of $G_0$.

A similar analysis can be performed for the right hand side of Eq.~(\ref{eq:nsbim}). For the first two terms on the right hand side, we can write with again two similar Taylor expansions: $\partial \phi(\boldsymbol x)/\partial n \approx \partial \phi (\boldsymbol x_0)/ \partial n + \nabla [\partial \phi/\partial n]\cdot (\boldsymbol x- \boldsymbol x_0)$ and $\partial f(\boldsymbol x)/ \partial n \approx \partial f(\boldsymbol x_0)/ \partial n +\nabla[\partial f/\partial n] \cdot (\boldsymbol x - \boldsymbol x_0)$, and then 
\begin{equation}
    \begin{aligned}
    &\frac{\partial \phi (\boldsymbol x)}{\partial n} G_k - \frac{\partial \phi(\boldsymbol x_0)}{\partial n}\frac{\partial f(\boldsymbol x)}{\partial n} G_0 
    \\
    \approx  &\frac{\partial \phi (\boldsymbol x_0)}{\partial n} [G_k - G_0] + \nabla \frac{\partial \phi}{\partial n} \cdot (\boldsymbol x - \boldsymbol x_0) G_k
    - \nabla \frac{\partial f}{\partial n} \cdot (\boldsymbol x - \boldsymbol x_0) G_0 
    \end{aligned}
\end{equation}
where we have used $\partial f(\boldsymbol x_0)/ \partial n = 1$. The term with $[G_k - G_0]$ is regular, and the terms with the gradients contain both $(\boldsymbol x - \boldsymbol x_0)$, which thus cancel out the $1/r$ singularity from $G_k$ and $G_0$, respectively. 

Since $f$ approaches zero as $\boldsymbol x \rightarrow \BS{x}_0$, it cancels out the third term with $\partial G_0/\partial n$ on the right hand side of Eq.~(\ref{eq:nsbim}).

In the Burton-Miller implementation of Eq.~(\ref{eq:d2Gkdndn0_d2G0dndn0}) (to be shown in Appendix~\ref{App:nsbim_BM}) it is claimed that
$$
\frac{\dff^2{G_{k}}}{\dff{n}\dff{n_0}}-\frac{\dff^2{G_{0}}}{\dff{n}\dff{n_0}} \rightarrow \frac{k^2}{2|\BS{x}-\BS{x}_0|}$$ 
when $\BS{x}\rightarrow\BS{x}_0$. The proof of this statement will be provided here. Due care has to be taken with the derivative $\partial / \partial n_0$, which now acts on $\boldsymbol x_0$ and not on $\boldsymbol x$, creating an additional minus sign. Applying this to Eq.~(\ref{eq:AppBdGdn}) we get:
\begin{equation}
\begin{aligned}
\frac{\partial^2 G_k}{\partial n_0 \partial n}&= \frac{e^{\rmi kr}}{r^3}\left\{-\boldsymbol n_0 \cdot \boldsymbol n (\rmi kr-1) - \frac{\boldsymbol n_0 \cdot (\boldsymbol x - \boldsymbol x_0) \boldsymbol n \cdot (\boldsymbol x - \boldsymbol x_0)}{r^2}(-k^2 r^2 -3 \rmi kr +3)\right\},
\\
\frac{\partial^2 G_0}{\partial n_0 \partial n}&=\boldsymbol n_0 \cdot \boldsymbol n\frac{1}{r^3} -3\frac{\boldsymbol n_0 \cdot (\boldsymbol x - \boldsymbol x_0) \boldsymbol n \cdot (\boldsymbol x - \boldsymbol x_0)}{r^5}.
\end{aligned}
\end{equation}
Using a Taylor expansion $e^{\rmi kr}\approx 1 + \rmi kr - k^2r^2/2$, we get:
\begin{equation}\label{eq:d2Gkdn0dn_d2G0dndn0Taylor}
\begin{aligned}
\frac{\partial^2 G_k}{\partial n_0 \partial n}-\frac{\partial^2 G_0}{\partial n_0 \partial n}= &\frac{1}{r^3}\boldsymbol n_0 \cdot \boldsymbol n \left(\frac{1}{2}k^2 r^2 + o(k^3r^3)\right)\\& - \frac{\boldsymbol n_0 \cdot (\boldsymbol x - \boldsymbol x_0) \boldsymbol n \cdot (\boldsymbol x - \boldsymbol x_0)}{r^5}\left(\frac{1}{2}k^2 r^2 +o(k^3r^3)\right).
\end{aligned}
\end{equation}
As explained before, when $\BS{x}\rightarrow\BS{x}_0$, term $\boldsymbol n_0 \cdot (\boldsymbol x - \boldsymbol x_0)$ and term $\boldsymbol n \cdot (\boldsymbol x - \boldsymbol x_0)$ both converge as $r^2$ and thus the last part of Eq.~(\ref{eq:d2Gkdn0dn_d2G0dndn0Taylor}) is non-singular. The first term on the right hand side gives the required $k^2/(2r)$ in Eq.~(\ref{eq:d2Gkdndn0_d2G0dndn0}), since $\boldsymbol n_0 \cdot \boldsymbol n \to 1$.

\section{Remarks on integrals at infinity and Sommerfeld condition}
\label{sec:app_inftyintegrals}
There are two remarks that can be made about the integrals at infinity that occur for external Helmholtz problems in the boundary element framework. The first one refers to the application of the Sommerfeld radiation condition~\cite{Sommerfeld1912, Schot1992} and applies to the conventional boundary element method as given in Eq.~(\ref{eq:BEM}) where the two integrals at infinity are:
\begin{equation}
    \begin{aligned}
    \int_\infty \phi \frac{\partial G_k}{\partial n} \rmd S - \int_\infty \frac{\partial\phi}{\partial n}  G_k \rmd S=\int_\infty \left[\phi \frac{\partial G_k}{\partial n}  -  \frac{\partial\phi}{\partial n}  G_k\right] \rmd S.
    \end{aligned}
\end{equation}
Suppose we draw a very large sphere with radius $R$ and consider this as `infinity' then $G_k = e^{\rmi kR}/R$ and $\partial G_k/\partial n=e^{\rmi kR} \boldsymbol n \cdot (\boldsymbol x - \boldsymbol x_0) (\rmi kR-1)/R^3$. Since $\boldsymbol n \cdot (\boldsymbol x - \boldsymbol x_0)\to R$ and $dS=4\pi R^2$, then:
\begin{equation}
    \begin{aligned}
    \int_\infty \left[\phi \frac{\partial G_k}{\partial n}  -  \frac{\partial\phi}{\partial n}  G_k\right] \rmd S
    \end{aligned} = \frac{e^{\rmi kR}}{R^3}\left[ R(\rmi kR-1)\phi - R^2 \frac{\partial \phi}{\partial n}\right]4\pi R^2.
\end{equation}
The term with $-1$ vanishes as $R\to \infty$ since $\phi$ decays at least as fast as $1/R$. Then:
\begin{equation}
    \begin{aligned}
    \int_\infty \left[\phi \frac{\partial G_k}{\partial n}  -  \frac{\partial\phi}{\partial n}  G_k\right] \rmd S
    \end{aligned} = e^{\rmi kR}R\left[ \rmi k\phi - \frac{\partial \phi}{\partial n}\right]4\pi
\end{equation}
which can be shown to be zero with the help of the Sommerfeld radiation condition~\cite{Sommerfeld1912, Schot1992} expressing the fact that only outgoing waves are allowed (and realizing that at $\infty$, $\partial \phi/\partial n = \partial \phi/ \partial r$) as
\begin{align}
    \lim_{r \to \infty}r\left[\frac{\dff \phi}{\dff r}-\rmi k \phi \right]  = 0. \label{eq:SommerfeldCon}
\end{align}
This is the formula as it appears in the original Sommerfeld article~\cite{Sommerfeld1912} (his Eq. (21)).

The second place where due care has to be taken with integrals at infinity is in Sec.~\ref{sec:BEMChapter}, where the desingularized boundary element method is described. In this case, there is actually a term that appears from one of the integrals at infinity. There are four integrals that now need to be considered in Eq.~(\ref{eq:nsbim}):
\begin{equation}
\begin{aligned}
    \int_{\infty}   \phi(\BS{x}_0) g(\BS{x}) \frac{\dff{G_{0}}}{\dff{n}}   \text{ d} S(\BS{x}) \qquad &; \qquad
    \int_{\infty}    \phi(\BS{x}_0) \frac{\dff{g(\BS{x})}}{\dff{n}} G_{0}  \text{ d} S(\BS{x}) \\
     \int_{\infty}    \frac{\dff{\phi(\BS{x}_0)}}{\dff{n}} \frac{\dff{f(\BS{x})}}{\dff{n}} G_{0}  \text{d} S(\BS{x})\qquad &; \qquad
    \int_{\infty}     \frac{\dff{\phi(\BS{x}_0)}}{\dff{n}} f(\BS{x}) \frac{\dff{G_{0}}}{\dff{n}}  \text{d} S(\BS{x}). 
\end{aligned}
\end{equation}
With the particular choice of $g=1$ from Eq.~(\ref{eq:nsbim_gf}), it can be seen that the second integral is zero immediately ($\partial g/ \partial n=0$). The first integral will become, with $\partial G_0/\partial n = (\boldsymbol x- \boldsymbol x_0) \cdot \boldsymbol n /r^3 \to R/R^3$ (as remarked earlier $\boldsymbol n \cdot (\boldsymbol x - \boldsymbol x_0)\to R$) and $\rmd S= 4 \pi R^2$, 
\begin{equation}
\begin{aligned}
    \int_{\infty}   \phi(\BS{x}_0) g(\BS{x}) \frac{\dff{G_{0}}}{\dff{n}}   \text{ d} S(\BS{x}) =\phi(\BS{x}_0) \frac{1}{R^3} R 4\pi R^2 = 4 \pi \phi(\BS{x}_0).
\end{aligned}
\end{equation}
This is the value mentioned in Sec.~\ref{sec:BEMChapter}. The other two remaining integrals will give zeroes since $f= \boldsymbol n(\boldsymbol x_0)\cdot (\boldsymbol x - \boldsymbol x_0) \to R \cos \theta$ and $\partial f/ \partial n = \boldsymbol n (\boldsymbol x_0) \cdot \boldsymbol n (\boldsymbol x) \to \cos \theta$, with $\theta$ the angle between the normal vector on the surface at $\boldsymbol x_0$ and the normal vector on the surface at infinity. This will cause the contributions of the integrals at infinity to cancel for these two last integrals. Thus: 
\begin{equation}
\begin{aligned}
     \int_{\infty}    \frac{\dff{\phi(\BS{x}_0)}}{\dff{n}} \frac{\dff{f(\BS{x})}}{\dff{n}} G_{0}  \text{d} S(\BS{x})=0
%\end{aligned}
%\end{equation}
\quad \text{and} \quad
%\begin{equation}
%\begin{aligned}
    \int_{\infty}     \frac{\dff{\phi(\BS{x}_0)}}{\dff{n}} f(\BS{x}) \frac{\dff{G_{0}}}{\dff{n}}  \text{d} S(\BS{x})=0. 
\end{aligned}
\end{equation}

\section{Derivation for the dielectric case}\label{App:Dielectric}

We need to expand $\boldsymbol n \cdot \partial \boldsymbol E/ \partial t_1$ and $\boldsymbol n \cdot \partial \boldsymbol E/ \partial t_2$ in Eq.~(\ref{eq:diel_Ht1_Ht2}). First we write again $\boldsymbol E = E_n \boldsymbol n + E_{t1} \boldsymbol t_1 + E_{t2} \boldsymbol t_2$. Since $\partial \boldsymbol n/ \partial t_1 = -\kappa_1 \boldsymbol t_1$ and $\partial \boldsymbol t_1/ \partial t_1 = \kappa_1 \boldsymbol n$ and $\partial \boldsymbol t_2/ \partial t_1 = 0$, we obtain: $\boldsymbol n \cdot \partial \boldsymbol E/ \partial t_1 = \boldsymbol n \cdot \boldsymbol n \partial E_n/\partial t_1 -\kappa_1 \boldsymbol n \cdot \boldsymbol t_1 E_n + \boldsymbol n \cdot \boldsymbol t_1 \partial E_{t1}/\partial t_1 + \kappa_1 \boldsymbol n \cdot \boldsymbol n E_{t1} + \boldsymbol n \cdot \boldsymbol t_2 \partial E_{t2}/ \partial t_1 +0 = \partial E_n/ \partial t_1 + \kappa_1 E_{t1}$.  We can then write:

\begin{equation}
\begin{aligned} \label{eq:diel_ndE/dt1}
    \boldsymbol n \cdot \frac{\partial \boldsymbol E}{\partial t_1} = \frac{\partial E_n}{\partial t_1}
 + \kappa_1 E_{t1}\quad \text{and similarly} \quad
 \boldsymbol n \cdot \frac{\partial \boldsymbol E}{\partial t_2} = \frac{\partial E_n}{\partial t_2}
 + \kappa_2 E_{t2}
 \end{aligned}
\end{equation}
where $\kappa_1$ is the curvature along the $\BS{t}_1$ direction and $\kappa_2$ the curvature along the $\BS{t}_2$ direction. We assumed that $\mu_{\rmin}= \mu_{\rmout}$, in Eq.~(\ref{eq:diel_Ht1_Ht2}) and using the tangential boundary conditions for the tangential magnetic field, Eq.~\ref{eq:diel_BC_Et1}, results in $\boldsymbol t_1 \cdot \frac{\partial \boldsymbol E^{\rmtr}}{\partial n}- \boldsymbol n \cdot \frac{\partial \boldsymbol E^{\rmtr}}{\partial t_1} = \boldsymbol t_1 \cdot \frac{\partial \boldsymbol E^{\rminc} +\boldsymbol E^{\rmsc}}{\partial n}- \boldsymbol n \cdot \frac{\partial \boldsymbol E^{\rminc}+\boldsymbol E^{\rmsc}}{\partial t_1}$. Thus Eq.~(\ref{eq:diel_ndE/dt1}) becomes
% %
% \begin{equation} \label{eq:diel_t1_dE/dn}
%     \boldsymbol t_1 \cdot \frac{\partial \boldsymbol E^{\rmtr}}{\partial n} - \frac{\partial E_n^{\rmtr}} {\partial t_1} - \kappa_1 E_{t1}^{\rmtr}  = \boldsymbol t_1 \cdot \frac{\partial \boldsymbol E^{\rmsc} -\boldsymbol E^{\rminc}}{\partial n} + \frac{\partial E_n^{\rmsc} + E_n^{\rminc}}{\partial t_1} - \kappa_1 (E_{t1}^{\rmsc} + E_{t1}^{\rminc}).
% \end{equation}
% %
%
\begin{equation}
\begin{aligned}
\label{eq:diel_t1_dE/dn}
    \boldsymbol t_1 \cdot \frac{\partial \boldsymbol E^{\rmtr}}{\partial n} - \frac{\partial E_n^{\rmtr}} {\partial t_1} - \kappa_1 E_{t1}^{\rmtr}  = \boldsymbol t_1 \cdot \frac{\partial (\boldsymbol E^{\rmsc} +\boldsymbol E^{\rminc})}{\partial n} - \frac{\partial (E_n^{\rmsc} + E_n^{\rminc})}{\partial t_1} - \kappa_1 (E_{t1}^{\rmsc} + E_{t1}^{\rminc})\\
    \boldsymbol t_2 \cdot \frac{\partial \boldsymbol E^{\rmtr}}{\partial n} - \frac{\partial E_n^{\rmtr}} {\partial t_2} - \kappa_2 E_{t2}^{\rmtr}  = \boldsymbol t_2 \cdot \frac{\partial (\boldsymbol E^{\rmsc} +\boldsymbol E^{\rminc})}{\partial n} - \frac{\partial (E_n^{\rmsc} + E_n^{\rminc})}{\partial t_2} - \kappa_2 (E_{t2}^{\rmsc} + E_{t2}^{\rminc}).
    \end{aligned}
\end{equation}
Due to the tangential continuity of the electric field in Eq.~(\ref{eq:diel_BC_Et1}), the terms with $\kappa_1$ and $\kappa_2$ disappear in the above equation. Further using the normal conditions on the electric field in Eq.~(\ref{eq:diel_BC_En}) to replace $E_n^{\rmtr}=(\epsilon_\rmout/\epsilon_\rmin)(E_n^\rmsc + E_n^\rminc)$ into Eq.~(\ref{eq:diel_t1_dE/dn}) gives:
\begin{equation}
\begin{aligned} \label{eq:diel_tdEtr/dn}
    \boldsymbol t_1 \cdot \frac{\partial \boldsymbol E^{\rmtr}}{\partial n}  = \boldsymbol t_1 \cdot \frac{\partial (\boldsymbol E^{\rmsc} +\boldsymbol E^{\rminc})}{\partial n} + \left[\frac{\epsilon_{\rmout}}{\epsilon_{\rmin}} -1\right] \frac{\partial (E_n^{\rmsc} + E_n^{\rminc})}{\partial t_1},\\
     \boldsymbol t_2 \cdot \frac{\partial \boldsymbol E^{\rmtr}}{\partial n}  = \boldsymbol t_2 \cdot \frac{\partial ( \boldsymbol E^{\rmsc} +\boldsymbol E^{\rminc})}{\partial n} + \left[\frac{\epsilon_{\rmout}}{\epsilon_{\rmin}} -1\right] \frac{\partial (E_n^{\rmsc} + E_n^{\rminc})}{\partial t_2}.
    \end{aligned}
\end{equation}
Thus the two tangential components of the vector $\partial \boldsymbol E^{\rmtr}/{\partial n}$ are now known in terms of scattered and incoming electric fields. Applying the divergence free condition of the scattered and total external field using Eq.~(\ref{eq:DivE_Surface}) gives:
\begin{equation}
\begin{aligned}
    &\boldsymbol n \cdot \frac{\partial \boldsymbol E^{\rmtr}}{\partial n} - \kappa E_n^{\rmtr} + \frac{\partial E_{t1}^{\rmtr}}{\partial t_1} + \frac{\partial E_{t2}^{\rmtr}}{\partial t_2} =0,\\
    &\boldsymbol n \cdot \frac{\partial (\boldsymbol E^{\rmsc}+\boldsymbol E^{\rminc})}{\partial n} - \kappa (E_n^{\rmsc}+E_n^{\rminc}) + \frac{\partial (E_{t1}^{\rmsc} +E_{t1}^{\rminc})}{\partial t_1} + \frac{\partial (E_{t2}^{\rmsc} +E_{t2}^{\rminc})}{\partial t_2}=0.
    \end{aligned}
\end{equation}
Subtracting these two equations and using the tangential boundary conditions in Eqs.~(\ref{eq:diel_BC_Et1})  again, will eliminate the terms with $\partial /\partial t_1$ and $\partial /\partial t_2$. Also using the normal conditions on the electric field to replace $E_n^{\rmtr}$:
\begin{equation} \label{eq:diel_ndEtr/dn}
    \boldsymbol n \cdot \frac{\partial \boldsymbol E^{\rmtr}}{\partial n}   =\boldsymbol n \cdot \frac{\partial \boldsymbol E^{\rmsc}}{\partial n} + \boldsymbol n \cdot \frac{\partial \boldsymbol E^{\rminc}}{\partial n}  + \kappa \left[\frac{\epsilon_{\rmout}}{\epsilon_{\rmin}}-1\right](E_n^{\rmsc}+E_n^{\rminc}).
\end{equation}
Now we have all the components of the normal derivative of the transmitted electric field, both in Cartesian as in tangential and normal decomposition. Let us write
\begin{equation}
    \begin{aligned}
    \frac{\partial \boldsymbol E^{\rmtr}}{\partial n} &= \boldsymbol e_x \frac{\partial E_x^{\rmtr}}{\partial n}
    +\boldsymbol e_y \frac{\partial E_y^{\rmtr}}{\partial n}
    +\boldsymbol e_z \frac{\partial E_z^{\rmtr}}{\partial n}
    \\
    &=\boldsymbol n \left(\boldsymbol n \cdot \frac{\partial \boldsymbol E^{\rmtr}}{\partial n}\right)
    +\boldsymbol t_1 \left(\boldsymbol t_1 \cdot \frac{\partial \boldsymbol E^{\rmtr}}{\partial n}\right)
    +\boldsymbol t_2 \left(\boldsymbol t_2 \cdot \frac{\partial \boldsymbol E^{\rmtr}}{\partial n}\right).
    \end{aligned}
\end{equation}
To reconstruct the Cartesian components, say the $x$-component of $\frac{\partial \boldsymbol E^{\rmtr}}{\partial n}$, we can use:
\begin{equation} \label{eq:diel_dEtrx/dn}
    \begin{aligned}
     \frac{\partial  E^{\rmtr}_x}{\partial  n}=\boldsymbol e_x \cdot \frac{\partial \boldsymbol E^{\rmtr}}{\partial n} = n_x \left(\boldsymbol n \cdot \frac{\partial \boldsymbol E^{\rmtr}}{\partial n}\right)
    + t_{1x} \left(\boldsymbol t_1 \cdot \frac{\partial \boldsymbol E^{\rmtr}}{\partial n}\right)
    + t_{2x} \left(\boldsymbol t_2 \cdot \frac{\partial \boldsymbol E^{\rmtr}}{\partial n}\right)
    \end{aligned}
\end{equation}
again with $n_x=(\boldsymbol n \cdot \boldsymbol e_x)$, $t_{1x}=(\boldsymbol t_1 \cdot \boldsymbol e_x)$ and $t_{2x}=(\boldsymbol t_2 \cdot \boldsymbol e_x)$. If we replace the terms in brackets in Eq.~(\ref{eq:diel_dEtrx/dn}) with Eq.~(\ref{eq:diel_ndEtr/dn}) and Eq.~(\ref{eq:diel_tdEtr/dn}) then:
\begin{equation} 
\label{eq:diel_dEtrx/dn2}
    \begin{aligned}
     \frac{\partial  E^{\rmtr}_x}{\partial  n}\; =\quad   n_x &\left\{
     \boldsymbol n \cdot \frac{\partial \boldsymbol E^{\rmsc}}{\partial n}
     +
      \boldsymbol n \cdot \frac{\partial \boldsymbol E^{\rminc}}{\partial n}
      +
      \kappa \left[\frac{\epsilon_{\rmout}}{\epsilon_{\rmin}}-1 \right](E_n^{\rmsc}+E_n^{\rminc})
     \right\}
     \\
      + t_{1x} &\left\{
     \boldsymbol t_1 \cdot \frac{\partial \boldsymbol E^{\rmsc}}{\partial n}
     +
      \boldsymbol t_1 \cdot \frac{\partial \boldsymbol E^{\rminc}}{\partial n}
      +
      \left[\frac{\epsilon_{\rmout}}{\epsilon_{\rmin}}-1 \right]\left(\frac{\partial E_n^{\rmsc}}{\partial t_1} + \frac{\partial E_n^{\rminc}}{\partial t_1}\right)
     \right\}
     \\
     + t_{2x} & \left\{
     \boldsymbol t_2 \cdot \frac{\partial \boldsymbol E^{\rmsc}}{\partial n}
     +
      \boldsymbol t_2 \cdot \frac{\partial \boldsymbol E^{\rminc}}{\partial n}
      +
      \left[\frac{\epsilon_{\rmout}}{\epsilon_{\rmin}}-1 \right]\left(\frac{\partial E_n^{\rmsc}}{\partial t_2} + \frac{\partial E_n^{\rminc}}{\partial t_2}\right)
     \right\}
     \\
     = \quad   n_x & \left(\boldsymbol n \cdot \frac{\partial \boldsymbol E^{\rmsc}}{\partial n}\right)
    + t_{1x} \left(\boldsymbol t_1 \cdot \frac{\partial \boldsymbol E^{\rmsc}}{\partial n}\right)
    + t_{2x} \left(\boldsymbol t_2 \cdot \frac{\partial \boldsymbol E^{\rmsc}}{\partial n}\right)
     + \frac{\partial E_x^{\rminc}}{\partial n}
     \\
     +\left[\frac{\epsilon_{\rmout}}{\epsilon_{\rmin}}-1 \right]& \left\{ 
     \kappa n_x E_n^{\rmsc} + t_{1x} \frac{\partial E_n^{\rmsc}}{\partial t_1}
     + t_{2x} \frac{\partial E_n^{\rmsc}}{\partial t_2}
     + \kappa n_x E_n^{\rminc} 
     + t_{1x} \frac{\partial E_n^{\rminc}}{\partial t_1}
     + t_{2x} \frac{\partial E_n^{\rminc}}{\partial t_2}
     \right\}.
\end{aligned}
\end{equation}
We express the transmitted field in the $x$-direction $E_x^{\rmtr}$ into normal and tangential components of the scattered and incoming field with the boundary conditions as:
\begin{equation} \label{eq:diel_Extr}
    E_x^{\rmtr} = \frac{\epsilon_{\rmout}}{\epsilon_{\rmin}} n_x E_n^{\rmsc} 
    +t_{1x} E_{t1}^{\rmsc}
    +t_{2x} E_{t2}^{\rmsc}
    +\frac{\epsilon_{\rmout}}{\epsilon_{\rmin}} n_x E_n^{\rminc} 
    +t_{1x} E_{t1}^{\rminc}
    +t_{2x} E_{t2}^{\rminc}.
\end{equation}
Finally we can perform the boundary element framework on the internal domain as ${\cal{H}}_{\rmin} \underline{E_x^{\rmtr}} ={\cal{G}}_{\rmin} \underline{\partial E_x^{\rmtr}/\partial n}$ and replace $E_x^{\rmtr}$ and its normal derivative by Eq.~(\ref{eq:diel_dEtrx/dn2}) and Eq.~(\ref{eq:diel_Extr}). Doing the same procedure for the $y$ and $z$ components and the external domain, a $6N \times 6N$ matrix will appear expressed in scattered electric field quantities: 

\begin{equation}\label{eq:big_dielectric}
    \begin{bmatrix}
      {\cal{H}} n_{x} & {\cal{H}} t_{1x} & {\cal{H}} t_{2x}
      & - {\cal{G}} n_x & - {\cal{G}} t_{1x}
      & - {\cal{G}} t_{2x} \\
      {\cal{H}} n_y & {\cal{H}} t_{1y}& {\cal{H}} t_{2y}
      & - {\cal{G}} n_y & - {\cal{G}} t_{1y}
      & - {\cal{G}} t_{2y}\\
      {\cal{H}} n_z & {\cal{H}} t_{1z} & {\cal{H}} t_{2z}
      & - {\cal{G}} n_z& - {\cal{G}} t_{1z}
      & - {\cal{G}} t_{2z}\\
       \bar{{\cal{H}}}^{n_{x}}_{\rmin}  & {\cal{H}}_{\rmin} t_{1x}
      &   {\cal{H}}_{\rmin} t_{2x} & - {\cal{G}}_{\rmin} n_x
      & - {\cal{G}}_{\rmin} t_{1x} & -  {\cal{G}}_{\rmin} t_{2x} \\
      \bar{{\cal{H}}}^{n_{y}}_{\rmin}  & {\cal{H}}_{\rmin} t_{1y}
      &   {\cal{H}}_{\rmin} t_{2y} &  - {\cal{G}}_{\rmin} n_y
      & - {\cal{G}}_{\rmin} t_{1y}& - {\cal{G}}_{\rmin} t_{2y}  \\ 
      \bar{{\cal{H}}}^{n_{z}}_{\rmin}  & {\cal{H}}_{\rmin} t_{1z}
      &   {\cal{H}}_{\rmin} t_{2z} &  - {\cal{G}}_{\rmin} n_z
      & - {\cal{G}}_{\rmin} t_{1z}& - {\cal{G}}_{\rmin} t_{2z}
    \end{bmatrix}
    \begin{bmatrix}
      E_n^{\rmsc} \\ E_{t_{1}}^{\rmsc} \\ E_{t_{2}}^{\rmsc} \\
      \boldsymbol{n} \cdot \frac{\partial \boldsymbol{E}^{\rmsc}}{\partial n}   \\
      \boldsymbol{t}_{1} \cdot \frac{\partial \boldsymbol{E}^{\rmsc}}{\partial n}  \\
      \boldsymbol{t}_{2} \cdot \frac{\partial \boldsymbol{E}^{\rmsc}}{\partial n}
    \end{bmatrix}
    =
    \begin{bmatrix}
      0 \\ 0 \\ 0 \\
      {\cal{B}}_{x} \\ {\cal{B}}_{y} \\ {\cal{B}}_{z}
    \end{bmatrix} 
  \end{equation}
  with $\alpha = x,\, y,$ or $z$ and
\begin{equation}
  \begin{aligned}
    \bar{{\cal{H}}}^{n_\alpha}_{\rmin} &= \frac{\epsilon_{\rmout}}{\epsilon_{\rmin}}  {\cal{H}}_{\rmin} n_{\alpha} - \left[\frac{\epsilon_{\rmout}}{\epsilon_{\rmin}}-1\right] {\cal{G}}_{\rmin} \left\{ \kappa 
    n_{\alpha} 
    + t_{1\alpha}
    \frac{\partial }{\partial t_1}
    + t_{2\alpha} 
    \frac{\partial}{\partial t_2} \right\} , \label{eq:matrix_partial_derivative} 
    \\
    {\cal{B}}_{\alpha} =& -\bar{{\cal{H}}}^{n_\alpha}_{\rmin} E_n^{\rminc} - {\cal{H}}_{\rmin} \left[t_{1\alpha} E^{\rminc}_{t_{1}} + t_{2\alpha}E^{\rminc}_{t_{2}}\right] \nonumber
  + {\cal{G}}_{\rmin} \frac{\partial E_\alpha^{\rminc}}{\partial n} .
  \end{aligned}
\end{equation}
The first row of the matrix system just states ${\cal{H}} \underline {E_x^{\rmsc}} = {\cal{G}} \underline {\partial E_x^{\rmsc}/ \partial n}$ (and similar for rows 2 and 3 for $y$ and $z$), but expressed in terms of normal and tangential components.

\section{Desingularization of the Burton-Miller method}\label{App:nsbim_BM}
We demonstrate here a Burton-Miller boundary integral method for the Helmholtz equation without any singularities. To the best of our knowledge a full desingularization of the Burton-Miller method has not appeared in literature yet, although hypersingular integrals have been studied extensively in the past~\cite{Chen1999}. 

We can perform the operator $\dff{(\cdot)}/\dff{n_0} \equiv \BS{n}(\BS{x}_{0})\BS{\cdot} \nabla_{\BS{x}_0}(\cdot)$. on Eq.~(\ref{eq:BEM})  with $\nabla_{\BS{x}_0}$ the gradient with respect to $\BS{x}_0$ to get: 
\begin{align}
   c(\BS{x}_0)\frac{\dff{\phi}(\BS{x}_{0})}{\dff{n}} + \int_{S} \phi(\BS{x}) \frac{\dff^2{G_{k}}}{\dff{n}\dff{n_0}} \text{ d} S(\BS{x}) = \int_{S}  \frac{\dff{\phi(\BS{x})}}{\dff{n}} \frac{\dff{G_{k}}}{\dff{n_0}} \text{ d} S(\BS{x}). \label{eq:BMcbim}
\end{align}
The term on the right hand side of Eq.~(\ref{eq:BMcbim}) is hyper-singular as $\BS{x}\rightarrow\BS{x}_{0}$. Performing the same operation on Eq.~(\ref{eq:cbimLapl}), we have
\begin{align}
    c(\BS{x}_0)\frac{ \dff{\psi(\BS{x}_{0}) } } {\dff{n}}  + \int_{S} \psi(\BS{x}) \frac{\dff^2{G_{0}}}{\dff{n}\dff{n_0}} \text{ d} S(\BS{x}) 
    = \int_{S}  \frac{ \dff{\psi(\BS{x})} } {\dff{n}} \frac{\dff{G_{0}}}{\dff{n_0}} \text{ d} S(\BS{x}). \label{eq:BMcbimLapl}
\end{align}
If we choose $\psi(\bs{x})=\phi(\bs{x}_0)$ in the domain which indicates $\nabla \psi(\bs{x}) = 0$, introduce it into Eq.~(\ref{eq:BMcbimLapl}), and subtract the resulting equation from Eq.~(\ref{eq:BMcbim}) we have
\begin{align}
    & c(\BS{x}_0)\frac{\dff{\phi}(\BS{x}_{0})}{\dff{n}}  + \int_{S} \phi(\BS{x}) \left[\frac{\dff^2{G_{k}}}{\dff{n}\dff{n_0}} - \frac{\dff^2{G_{0}}}{\dff{n}\dff{n_0}} \right] \text{ d} S(\BS{x}) = \int_{S}  \frac{\dff{\phi(\BS{x})}}{\dff{n}} \frac{\dff{G_{k}}}{\dff{n_0}} \text{ d} S(\BS{x}). \label{eq:BMcbimWeak}
\end{align}
The integrands in Eq.~(\ref{eq:BMcbimWeak}) are all weakly singular since
%:\\ \textcolor{red}{$\phi(\BS{x}) \frac{\dff^2{G_{k}}}{\dff{n}\dff{n_0}} - \phi(\BS{x}_{0}) \frac{\dff^2{G_{0}}}{\dff{n}\dff{n_0}}=[\phi(\BS{x})-\phi(\BS{x}_0)] \frac{\dff^2{G_{k}}}{\dff{n}\dff{n_0}} + \phi(\BS{x}_{0})(\frac{\dff^2{G_{k}}}{\dff{n}\dff{n_0}}-  \frac{\dff^2{G_{0}}}{\dff{n}\dff{n_0}})$} and 
(see Appendix~\ref{sec:app_notesDesing})
\begin{align} \label{eq:d2Gkdndn0_d2G0dndn0}
    \lim_{\BS{x} \to \BS{x}_0} \left( \frac{\dff^2{G_{k}}}{\dff{n}\dff{n_0}}-\frac{\dff^2{G_{0}}}{\dff{n}\dff{n_0}} \right) = \frac{k^2}{2|\BS{x}-\BS{x}_0|}.
\end{align}
To analytically remove the remaining singularities and the terms associated with the solid angle at $\BS{x}_0$ in Eq.~(\ref{eq:BMcbimWeak}), we follow the same procedure illustrated in Sec.~\ref{sec:BEMChapter} to set up the following two functions that satisfy the Laplace equation in the domain where $\phi$ is valid:
\begin{subequations}\label{eq:nsbim_BMpsi}
\begin{eqnarray} 
\psi_1({\BS{x}}) &=& \frac{k^2}{2} \phi(\BS{x}_0)  \, \left[\BS{n}(\BS{x}_0)\BS{\cdot}(\BS{x}-\BS{x}_0)\right], \\
\psi_2({\BS{x}}) &=& \frac{\dff{\phi(\BS{x}_0)}}{\dff{n}}.
\end{eqnarray}
\end{subequations}
Subtracting the conventional boundary integral equations corresponding to Eq.~(\ref{eq:nsbim_BMpsi}) from Eq.~(\ref{eq:BMcbimWeak}), and considering the integrals over the surface at $\infty$ for external problems, we get
\begin{align}
    &  \int_{S} \bigg\{ \phi(\BS{x}) \left[\frac{\dff^2{G_{k}}}{\dff{n}\dff{n_0}} - \frac{\dff^2{G_{0}}}{\dff{n}\dff{n_0}} \right] - \frac{k^2}{2} \phi(\BS{x}_{0}) [\BS{n}(\BS{x})\BS{\cdot} \BS{n}(\BS{x}_0)] G_{0} \bigg\} \text{d} S(\BS{x}) \nonumber \\
   = & \int_{S}  \left[ \frac{\dff{\phi(\BS{x})}}{\dff{n}} \frac{\dff{G_{k}}}{\dff{n_0}} + \frac{\dff{\phi(\BS{x}_0)}}{\dff{n}} \frac{\dff{G_{0}}}{\dff{n}}\right] \text{ d} S(\BS{x}) - \int_{S} \frac{k^2}{2} \phi(\BS{x}_0)  \, \left[\BS{n}(\BS{x}_0)\BS{\cdot}(\BS{x}-\BS{x}_0)\right]\frac{\dff{ G_{0}}}{\dff{n}} \text{ d} S(\BS{x}) -4\pi \frac{\dff{\phi(\BS{x}_0)}}{\dff{n}}.\label{eq:nsbimBMhyper}
\end{align}
The integrands in Eq.~(\ref{eq:nsbimBMhyper}) are now all regular and this equation is now fully desingularized.

\begin{figure}[t] 
\centering
\subfloat[Non-singular BEM of Eq.~(\ref{eq:nsbim})]{\includegraphics[width=0.45\textwidth]{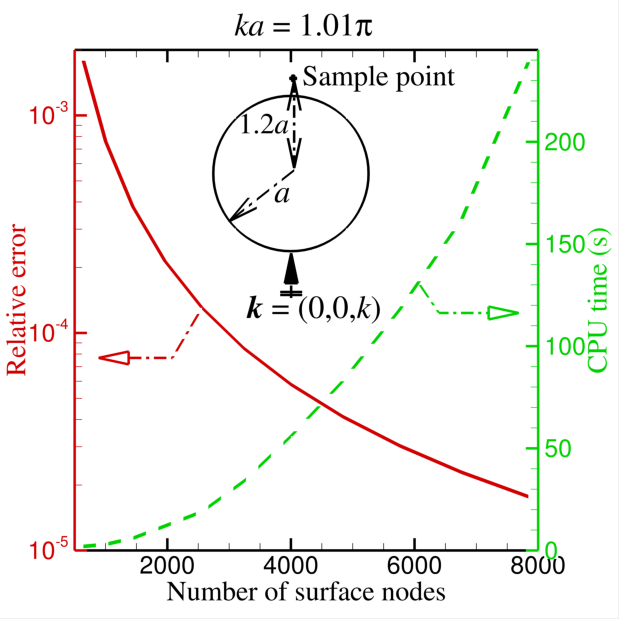}}  \qquad
\subfloat[Non-singular Burton-Miller BEM of Eq. ~(\ref{eq:nsbimBM_a})]{\includegraphics[width=0.45\textwidth]{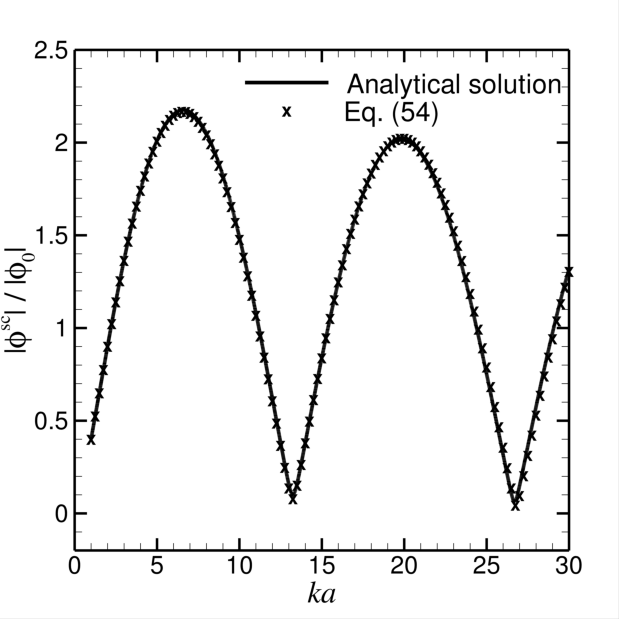}}
\caption{(a) Convergence and accuracy study for a standard scattering hard sound sphere at $ka=1.01\pi$ (very close to the internal resonance frequency of $ka=\pi$). The incoming wave travels from bottom to top and the observation point is taken at $r=1.2a$ above the sphere. The relative error between the numerical results by Eq.~(\ref{eq:nsbim}) and the analytical solution~\cite{Morse1991} reduces from $10^{-3}$ for 1000 nodes to about $10^{-5}$ for 8000 nodes for the scattered field. The computational time increases from a few seconds to 250 s. (b) Comparison between the analytical solution and the result with the non-singular Burton-Miller boundary element method of Eq.~(\ref{eq:nsbimBM_a}) for the case sketched in the inset of (a) when $ka$ is swept from 1 to 30 with a step of $\Delta ka =0.01$ (note only 1 in every 25 points is plotted). The sphere surface is represented by 7842 nodes connected by 3920 quadratic triangular elements. Good agreement is found and no fictitious frequencies appear. \label{fig:validation}}
\end{figure}

The Burton-Miller idea is to eliminate the spurious solutions at fictitious frequencies by combining Eq.~(\ref{eq:nsbimBMhyper}) and Eq.~(\ref{eq:nsbim}) (multiplied by an imaginary parameter $\rmi \beta$) since the fictitious frequencies in Eq.~(\ref{eq:nsbimBMhyper}) and Eq.~(\ref{eq:nsbim}) are always different from each other. To balance Eq.~(\ref{eq:nsbimBMhyper}) and Eq.~(\ref{eq:nsbim}) in terms of physical dimensions (dimensional homogeneous in length), the parameter $\beta$ should be related to a characteristic length of the problem under consideration. In our simulations, shown in Fig.~\ref{fig:validation}b, we used $\text{min}[0.5 a, 1/k]$. The obvious choices are the size of the scattering object or the inverse of the wave number $k$. If we introduce Eq.~(\ref{eq:nsbim_gf}) into Eq.~(\ref{eq:nsbim}), and combine the resulting formulation with Eq.~(\ref{eq:nsbimBMhyper}) according to the Burton-Miller idea, we obtain the non-singular Burton-Miller boundary integral method for the Helmholtz equation as
\begin{align}\label{eq:nsbimBM_a}
    &\quad \text{ }  4\pi \phi(\BS{x}_0) +  \int_{S} \left[\phi(\BS{x}) \frac{\dff{G_{k}}}{\dff{n}} - \phi(\BS{x}_0) \frac{\dff{G_{0}}}{\dff{n}}  \right] \text{ d} S(\BS{x}) \nonumber \\
    &\text{ }+\rmi \beta \int_{S} \bigg\{ \phi(\BS{x}) \left[\frac{\dff^2{G_{k}}}{\dff{n}\dff{n_0}} - \frac{\dff^2{G_{0}}}{\dff{n}\dff{n_0}} \right] - \frac{k^2}{2} \phi(\BS{x}_{0}) [\BS{n}(\BS{x})\BS{\cdot} \BS{n}(\BS{x}_0)] G_{0} \bigg\} \text{d} S(\BS{x}) \nonumber \\
    &= 
    \text{ } \int_{S}  \left[ \frac{\dff{\phi(\BS{x})}}{\dff{n}} G_{k} - \frac{\dff{\phi(\BS{x}_0)}}{\dff{n}} [\BS{n}(\BS{x}_0)\BS{\cdot} \BS{n}(\BS{x})] G_{0} \right] \text{d} S(\BS{x}) +\int_{S} \frac{\dff{\phi(\BS{x}_0)}}{\dff{n}} [\BS{n}(\BS{x}_0)\BS{\cdot} (\BS{x}-\BS{x}_0)] \frac{\dff{G_{0}}}{\dff{n}} \text{d} S(\BS{x}) \nonumber \\
    & \text{ } + \rmi \beta \int_{S}  \left[ \frac{\dff{\phi(\BS{x})}}{\dff{n}} \frac{\dff{G_{k}}}{\dff{n_0}} + \frac{\dff{\phi(\BS{x}_0)}}{\dff{n}} \frac{\dff{G_{0}}}{\dff{n}}\right] \text{ d} S(\BS{x})  - \rmi \beta \int_{S} \frac{k^2}{2} \phi(\BS{x}_0)  \, \left[\BS{n}(\BS{x}_0)\BS{\cdot}(\BS{x}-\BS{x}_0)\right]\frac{\dff{ G_{0}}}{\dff{n}} \text{ d} S(\BS{x}) \nonumber \\
    & \text{ } -4\pi \rmi \beta \frac{\dff{\phi(\BS{x}_0)}}{\dff{n}}.
\end{align}
It is worth noting that in Eq.~(\ref{eq:nsbimBM_a}), the integrals over the surface at $\infty$ for external problems have been included, resulting in the terms with $4\pi$. For internal problems, those terms disappear. Our goal to have a non-singular version of the Burton-Miller boundary integral method, as shown in Eq.~(\ref{eq:nsbimBM_a}) is now achieved and the fictitious frequencies are eliminated.

\bibliographystyle{ieeetr}
\bibliography{references}

\end{document}